\documentclass{amsart}

\newtheorem{theorem}{Theorem}[section]
\newtheorem{lemma}[theorem]{Lemma}

\theoremstyle{definition}

\newtheorem{example}[theorem]{Example}

\theoremstyle{remark}

\usepackage{url}

\begin{document}

\title{Landen Survey}

\author{Dante V. Manna}
\address{Department of Mathematics and Statistics, Dalhousie University, Halifax, Nova Scotia, Canada, B3H 3J5}
\email{dantemanna@gmail.com}

\author{Victor H. Moll}
\address{Department of Mathematics,
Tulane University, New Orleans, LA 70118}
\email{vhm@math.tulane.edu}


\date{\today}


\begin{abstract}
Landen transformations are maps on the coefficients of an integral that
preserve its value. We present a brief survey of their appearance in the 
literature.
\end{abstract}

\maketitle

\begin{center}
{\em To Henry, who provides inspiration, taste and friendship}
\end{center}

\newcommand{\nn}{\nonumber}
\newcommand{\ba}{\begin{eqnarray}}
\newcommand{\ea}{\end{eqnarray}}
\newcommand{\ift}{\int_{0}^{\infty}}
\newcommand{\ifft}{\int_{- \infty}^{\infty}}
\newcommand{\agm}{\text{AGM}}
\newcommand{\no}{\noindent}
\newcommand{\realpart}{\mathop{\rm Re}\nolimits}
\newcommand{\imagpart}{\mathop{\rm Im}\nolimits}

\newtheorem{Definition}{\bf Definition}[section]
\newtheorem{Thm}[Definition]{\bf Theorem} 
\newtheorem{Example}[Definition]{\bf Example} 
\newtheorem{Lem}[Definition]{\bf Lemma} 
\newtheorem{Note}[Definition]{\bf Note} 
\newtheorem{Cor}[Definition]{\bf Corollary} 
\newtheorem{Prop}[Definition]{\bf Proposition} 
\newtheorem{Problem}[Definition]{\bf Problem} 
\numberwithin{equation}{section}

\section{In the beginning there was Gauss} \label{sec-intro} 
\setcounter{equation}{0}

In the year 1985, one of us had the luxury of attending a graduate course 
on {\em Elliptic Functions} given by Henry McKean 
at the Courant Institute. Among 
the many beautiful results he described in his unique style, there was
a calculation of Gauss:  take two positive real numbers $a$ and $b$, with 
$a > b$,  and 
form a new pair by replacing $a$ with the arithmetic mean $(a+b)/2$
and $b$ with the geometric mean $\sqrt{ab}$. Then iterate:
\begin{equation}
a_{n+1} = \frac{a_{n}+b_{n}}{2}, \quad b_{n+1} = \sqrt{a_{n}b_{n}}
\label{agm-1}
\end{equation}
\noindent
starting with $a_{0}=a$ and $b_{0}=b$. Gauss \cite{gauss1} was 
interested in the initial conditions $a=1$ and $b= \sqrt{2}$.
The iteration generates a sequence of algebraic numbers which rapidly become
impossible to describe explicitly, for instance,
\begin{equation}
a_{3} = \frac{1}{2^{3}} \left( (1 + \sqrt[4]{2} )^{2} + 
2 \sqrt{2} \sqrt[8]{2} \, \sqrt{1+ \sqrt{2}} \right)
\end{equation}
\noindent
is a root of the polynomial 
\begin{eqnarray}
G(a) & = & 16777216a^8-16777216a^7+5242880a^6-10747904a^5 \nonumber \\
 & & +942080a^4-1896448a^3+4436a^2-59840a+1. \nonumber 
\end{eqnarray}
\noindent
The numerical behavior is surprising; $a_{6}$ and $b_{6}$ agree to 87 
digits. It is simple to check that 
\begin{equation}
\lim\limits_{n \to \infty} a_{n} =
\lim\limits_{n \to \infty} b_{n}. 
\end{equation}
\noindent
See (\ref{conv-elliptic}) for details. This common limit 
is called the {\em arithmetic-geometric mean} and is 
denoted by $\agm(a,b)$. It is the explicit dependence on the initial
condition that is hard to discover. 

Gauss computed some numerical values and observed that 
\begin{equation}
a_{11} \sim b_{11} \sim  1.198140235
\label{once}
\end{equation}
\noindent
and then he {\em recognized} the reciprocal of this number as a numerical 
approximation to the elliptic integral
\begin{equation}
I = \frac{2}{\pi} \int_{0}^{1} \frac{dt}{\sqrt{1-t^{4}}}. 
\end{equation}
\noindent
It is unclear to the authors how Gauss
recognized this number - he simply knew it.  (Stirling's tables may have been 
a help; \cite{borw1} contains a 
reproduction of the original notes and comments.) 
He was particularly interested in the evaluation of this definite integral
as it provides the length of a lemniscate.  
In his diary Gauss remarked,
{\em `This will surely open up a 
whole new field of analysis'} \cite{cox84, borwein1}. 

Gauss' procedure to find an analytic expression for $\agm(a,b)$ began with 
the elementary observation
\begin{equation}
\agm(a,b) = \agm \left( \frac{a+b}{2}, \sqrt{ab} \right)
\label{itera}
\end{equation}
\noindent
and the homogeneity condition
\begin{equation}
\agm(\lambda a, \lambda b)  = \lambda \agm(a,b)~.
\end{equation}
\noindent
He used (\ref{itera}) with $a= (1+ \sqrt{k})^{2}$ and 
$b= (1- \sqrt{k})^{2}$, with $0 < k < 1$, to produce
\begin{equation}
\agm ( 1+k + 2 \sqrt{k}, 1+k - 2 \sqrt{k} )
 = \agm (1+k,1-k). \nonumber 
\end{equation}
\noindent
He then used the homogeneity of $\agm \, $ to write
\begin{eqnarray}
\agm ( 1+k + 2 \sqrt{k}, 1+k - 2 \sqrt{k} )
               & = & \agm( (1+k)(1+k^{*}), (1+k)(1-k^{*})) \nonumber \\
               & = & (1+k) \agm( 1 + k^{*}, 1- k^{*} ), \nonumber 
\end{eqnarray}
\noindent
with 
\begin{equation}
k^{*}  = \frac{2 \sqrt{k}}{1+k}. 
\end{equation}
\noindent
This resulted in the functional equation
\begin{equation}
\agm ( 1 + k, 1-k) = (1+k) \, \agm ( 1 + k^{*}, 1 - k^{*}). \label{funct}
\end{equation}

In his analysis of (\ref{funct}), Gauss substituted the power series
\begin{equation}
\frac{1}{\agm(1+k,1-k)} = \sum_{n=0}^{\infty} a_{n}k^{2n} 
\end{equation}
\noindent
into  (\ref{funct}) and, solved an infinite system of
nonlinear equations,  to produce
\begin{equation}
a_{n} = 2^{-2n} \binom{2n}{n}^{2}.
\end{equation}
\noindent
Then he recognized the series as that of an elliptic integral to obtain 
\begin{equation}\label{agmlaw}
\frac{1}{\agm(1+k,1-k)} = \frac{2}{\pi} \int_{0}^{\pi/2} \frac{dx}
{\sqrt{1- k^{2}\sin^{2}x}}.
\end{equation}
\noindent
This is a remarkable tour de force.  \\

The function 
\begin{equation}
K(k) = \int_{0}^{\pi/2} \frac{dx}{\sqrt{1 - k^{2} \sin^{2} x} }
\end{equation}
\noindent 
is the {\em elliptic integral of the first kind}. It can also be written in 
the algebraic form
\begin{equation}
K(k) = \int_{0}^{1} \frac{dt}{\sqrt{(1-t^{2})(1-k^{2}t^{2})}}.
\end{equation}
\noindent
In this notation, (\ref{funct}) becomes 
\begin{equation}
K(k^{*}) = ( 1 + k) K(k).
\end{equation}

This is the {\em Landen transformation} for the complete elliptic 
integral.  John Landen \cite{landen1}, the namesake of the 
transformation, studied related integrals: for example,
\begin{equation} 
\kappa := \int_0^1 \frac{dx}{\sqrt{x^2(1-x^2)}}~. 
\end{equation} 
He derived identites such as 
\begin{equation} 
\kappa = \varepsilon \sqrt{\varepsilon^2-\pi}~,~\mbox{~where~~} 
\varepsilon := \int_0^{\pi/2} \sqrt{2- \sin^2 \theta}~ d\theta~, 
\end{equation} 
proven mainly by suitable changes of varibles in the integral 
for $\varepsilon$.  In \cite{watson1933} the reader will find a historical 
account of Landen's work, including the above identities.

The reader will find in \cite{borwein1} and 
\cite{mckmoll} proofs in a variety of styles. In trigonometric form, the 
Landen transformation states that 
\begin{equation}
G(a,b) = \int_{0}^{\pi/2} \frac{d \theta}{\sqrt{a^{2} \cos^{2} \theta 
+ b^{2} \sin^{2} \theta}} 
\label{elint-1}
\end{equation}
\noindent
is invariant under the change of parameters $(a,b) \mapsto 
\left( \tfrac{a+b}{2}, \sqrt{ab} \right)$. D. J. Newman \cite{newman1} 
presents a very clever proof: the change of variables $x = b \tan \theta$ 
yields 
\begin{equation}
G(a,b) = \frac{1}{2} \int_{-\infty}^{\infty} \frac{dx}
{\sqrt{ (a^{2}+x^{2})(b^{2}+x^{2})}}. 
\end{equation}
\noindent
Now let  $x \mapsto  x + \sqrt{x^{2} + ab}$  to complete the proof. Many of 
the above identities can now be searched for and proven on a 
computer \cite{borw1}.

\section{An interlude: the quartic integral} \label{sec-quartic} 
\setcounter{equation}{0}

The evaluation of definite integrals of rational functions is one of the
standard topics in Integral Calculus. Motivated by the lack of success of
symbolic languages, we began a systematic study of these 
integrals. {\it A posteriori}, one learns that even rational functions are 
easier to deal with. Thus we start with one having 
{\em a power of a quartic} in its denominator. The evaluation of the identity 
\begin{equation}
\ift \frac{dx}{(x^{4}+2ax^{2}+1)^{m+1} } 
= \frac{\pi}{2^{m+3/2} \, (a+1)^{m+1/2}} P_{m}(a),
\label{quartic}
\end{equation}
\noindent
where 
\begin{equation}
P_{m}(a) = \sum_{l=0}^{m} d_{l}(m) a^{l}
\end{equation}
\noindent
with
\begin{equation}
d_{l}(m) = 2^{-2m} \sum_{k=l}^{m} 2^{k} \binom{2m-2k}{m-k} \binom{m+k}{m} 
\binom{k}{l},
\label{dofl}
\end{equation}
\noindent
was first established in \cite{bomohyper}. 

A standard hypergeometric argument yields 
\begin{equation}
P_{m}(a) = P_{m}^{(\alpha, \beta)}(a), 
\end{equation}
\noindent
where 
\begin{equation}
P_{m}^{(\alpha,\beta)}(a) = \sum_{k=0}^{m} (-1)^{m-k} 
\binom{m+\beta}{m-k} \binom{m+k+\alpha + \beta}{k} 2^{-k} (a+1)^{k}
\end{equation}
\noindent
is the classical Jacobi polynomial; the parameters $\alpha$ and $\beta$
are given by 
$\alpha = m + \tfrac{1}{2}$ and $\beta = -m - \tfrac{1}{2}$. A general 
description of these functions and their properties are given 
in \cite{abramowitz1}. The twist here is that 
they depend on $m$, which means  most 
of the properties of $P_{m}$ had to be proven 
from scratch. For instance, $P_{m}$ satisfies the recurrence
\begin{eqnarray}
P_{m}(a) & = & \frac{(2m-3)(4m-3)a}{4m(m-1)(a-1)}P_{m-2}(a) - 
\frac{(4m-3)a(a+1)}{2m(m-1)(a-1)} P_{m-2}'(a) \nonumber \\
& + & \frac{4m(a^{2}-1)+1-2a^{2}}{2m(a-1)}P_{m}(a).  \nonumber 
\end{eqnarray}
\noindent
This {\em cannot} be obtained by replacing $\alpha = m + \tfrac{1}{2}$ and 
$\beta = -m - \tfrac{1}{2}$ in the standard recurrence for the Jacobi 
polynomials. The reader will find in \cite{tv1} several differnt proofs of 
(\ref{quartic}). \\

The polynomials $P_{m}(a)$ makes a surprising appearance in 
the expansion
\begin{equation}
\sqrt{a + \sqrt{1+c}} = \sqrt{a+1} \left[ 
1 - \sum_{k=1}^{\infty} \frac{(-1)^{k}}{k} \frac{P_{k-1}(a) \, c^{k} }
{2^{k+1} \, (a+1)^{k} } \right]
\end{equation}
\noindent
as described in \cite{bomoram}. The special case $a=1$ appears in 
\cite{bromwich}, page $191$, exercise $21$. Ramanujan \cite{berndt1} had a more 
general expression, but only for the case $c = a^{2}$:
\begin{equation}
(a  + \sqrt{1+a^{2}} )^{n} = 1+ na + \sum_{k=2}^{\infty} \frac{b_{k}(n)a^{k}}
{k!},
\end{equation}
\noindent
where, for $k \geq 2$, 
\begin{equation}
b_{k}(n) = \begin{cases}
            n^{2}(n^{2}-2^{2})(n^{2}-4^{2}) \cdots (n^{2} - (k-2)^{2}) 
\text{ if }k \text{ is even}, \\
            n(n^{2}-1^{2})(n^{2}-3^{2}) \cdots (n^{2} - (k-2)^{2}) 
\text{  if }k \text{ is odd}. 
\end{cases}
\end{equation}
\noindent
This result appears in Berndt \cite{berndt1} as Corollary 2 to Entry 14 
and is machine-checkable, as are many of the identities in this section.  \\

The coefficients $d_{l}(m)$ in (\ref{dofl}) have 
:many interesting properties: \\

\noindent
$\bullet$ They form a {\em unimodal sequence}: there exists an index 
$0 \leq m^{*} \leq m$ such that $d_{j}(m)$ increases up to $j=m^{*}$ and
decreases from then on. See \cite{bomouni1} for a proof of the more
general statement: {\em If} $P(x)$ {\em is 
a polynomial with nondecreasing, nonnegative
coefficients, then the coefficient sequence of} $P(x+1)$ 
{\em is unimodal}.  \\

\noindent
$\bullet$ They form a {\em log-concave sequence}: define the operator 
$\mathfrak{L}( \{ a_{k} \}) := \{ a_{k}^{2}-a_{k-1}a_{k+1} \}$ acting on 
sequences of positive real numbers. A sequence 
$\{ a_{k} \}$ is called log-concave if its image under $\mathfrak{L}$ is 
again a sequence of positive numbers; i.e. $a_{k}^{2} - a_{k-1}a_{k+1} 
\geq 0$. Note that this condition is satisfied if and only if the sequence 
$\{b_k := \log(a_k) \}$ is concave, hence the name.  We refer the reader 
to \cite{wilf1} for a detailed introduction.  The log-concavity of 
$d_{l}(m)$ was established in \cite{kauers-paule} using Computer 
Algebra techniques:  in particular, cylindrical algebraic decompositions 
as developed in \cite{caviness1} and \cite{collins2}.  \\

\noindent
$\bullet$ They produce interesting polynomials: in \cite{bomosha} 
one finds the representation
\begin{equation}
d_{l}(m) = \frac{A_{l,m}}{l! \, m! \, 2^{m+l}}, 
\end{equation}
\noindent
with
\begin{equation}
A_{l,m} = \alpha_{l}(m) \prod_{k=1}^{m} (4k-1) - 
           \beta_{l}(m) \prod_{k=1}^{m} (4k+1).
\end{equation}
\noindent
Here $\alpha_{l}$ and $\beta_{l}$ are polynomials in $m$ of degrees $l$ and 
$l-1$, respectively. For 
example, $\alpha_{1}(m) = 2m+1$ and $\beta_{1}(m) 
= 1$, so that the coefficient of the linear term of $P_{m}(a)$ is 
\begin{equation}
d_{1}(m) = \frac{1}{\, m! \, 2^{m+1}} 
\left( (2m+1)  \prod_{k=1}^{m} (4k-1) - 
            \prod_{k=1}^{m} (4k+1) \right). 
\end{equation}
\noindent
J. Little established in \cite{little} the remarkable fact that the 
polynomials $\alpha_{l}(m)$ and $\beta_{l}(m)$ have
all their roots on the vertical line $\realpart{m} = - \tfrac{1}{2}$.

When we showed this to Henry, he simply remarked: {\em the only thing you have 
to do now is to let $l \to \infty$ and get the Riemann hypothesis}. The proof
in \cite{little} consists in a study of the recurrence 
\begin{equation}
y_{l+1}(s) = 2sy_{l}(s) - \left( s^{2} -(2l-1)^{2} \right)y_{l-1}(s),
\end{equation}
\noindent
satisfied by $\alpha_{l}((s-1)/2)$ and 
$\beta_{l}((s-1)/2)$. There is no Number Theory in the proof, so 
it is not likely to connect to the Riemann zeta function 
$\zeta(s)$, but one never knows. 

The arithmetical properties of $A_{l,m}$ are beginning to be elucidated.
We have shown that their $2$-adic valuation satisfies 
\begin{equation}
\nu_{2}(A_{l,m}) = \nu_{2}( (m+1-l)_{2l}) + l,
\end{equation}
\noindent
where $(a)_{k} = a(a+1)(a+2) \cdots (a+k-1)$ is the Pochhammer symbol. This
expression allows for a combinatorial interpretation of the block 
structure of these valuations. See \cite{amm1} for details.

\section{The incipient rational landen transformation} \label{sec-ratlanden} 
\setcounter{equation}{0}

The clean analytic expression in (\ref{quartic}) is not expected to extend 
to rational
functions of higher order. In our analysis we distinguish according to the
domain of integration: the finite interval case,  mapped by a bilinear
transformation to $[0, \infty)$, and the whole line. In this section
we consider the definite integral,
\begin{equation}
U_{6}(a,b;c,d,e) = \ift \frac{cx^{4} + dx^{2} + e}{x^{6}+ax^{4}+bx^2+1} \, dx,
\label{usix}
\end{equation}
\noindent
as the simplest case on $[0, \infty)$.  The case of the real line is
considered below. The integrand is chosen to be 
even by necessity: {\em none of the techniques in 
this section work for the odd case}.
We normalize two of the coefficients in the denominator in order to reduce the 
number of parameters. The standard approach for the evaluation of 
(\ref{usix}) is to introduce the 
change of variables $x  = \tan \theta$. This leads to an 
intractable trigonometric integral. 

A different result is obtained if one 
first symmetrizes the denominator: we say that 
a polynomial of degree $d$ is {\em reciprocal} if $Q_{d}(1/x) = 
x^{-d}Q_{d}(x)$, that is, the sequence of its coefficients is a 
palindrome. Observe that if $Q_{d}$ is any polynomial of degree $d$, then 
\begin{equation}
T_{2d}(x) = x^{d}Q_{d}(x) Q_{d}(1/x)
\end{equation}
\noindent
is a reciprocal polynomial of degree $2d$.   For example, if 
\begin{equation}
Q_{6}(x) = x^{6} + ax^{4} + bx^{2} + 1. 
\end{equation}
\noindent
then 
\begin{eqnarray}
T_{12}(x) & = & x^{12} + (a+b)x^{10} + (a+b+ab)x^{8} + \nonumber \\
 & + & (2 + a^{2}+b^{2})x^{6} + 
(a+b+ab)x^{4} + (a+b)x^{2} + 1. 
\nonumber
\end{eqnarray}
\noindent
The numerator and denominator in the integrand of (\ref{usix}) are 
now scaled by $x^{6}Q_{6}(1/x)$ to produce a new integrand with reciprocal 
denominator: 
\begin{equation}
U_{6} = \ift \frac{S_{10}(x)}{T_{12}(x)} \, dx,
\end{equation}
\noindent
where we write
\begin{equation}
S_{10}(x) = \sum_{j=0}^{5} s_{j} x^{2j} \text{ and }
T_{12}(x) = \sum_{j=0}^{6} t_{j} x^{2j}. \nonumber 
\end{equation}
\noindent
The change of variables $x= \tan \theta$ now yields
\begin{equation}
U_{6} = \int_{0}^{\pi/2} \frac{S_{10}(\tan \theta) \, \cos^{10}(\theta)}
                              {T_{12}(\tan \theta) \, \cos^{12}(\theta)} \, 
d \theta. 
\end{equation}
\noindent
Now let $w = \cos 2 \theta$ and use $\sin^{2}\theta  = \frac{1}{2}(1-w)$
and $\cos^{2}\theta  = \frac{1}{2}(1+w)$ to check that the numerator 
and denominator of the new integrand,
\begin{equation}
S_{10}(\tan \theta) \, \cos^{10} \theta  = \sum_{j=0}^{5} s_{j} 
\sin^{2j}\theta \, \cos^{10-2j}\theta
\label{polyinw}
\end{equation}
\noindent
and 
\begin{eqnarray}
T_{12}(\tan \theta) \, \cos^{12} \theta  & = & \sum_{j=0}^{6} t_{j} 
\sin^{2j}\theta \, \cos^{12-2j}\theta  \nonumber \\
& = & 2^{-6} \sum_{j=0}^{6} t_{j} (1 - w)^{j} (1+w)^{6-j},
\nonumber
\end{eqnarray}
\noindent
are both polynomials in $w$. The mirror symmetry of $T_{12}$, reflected in 
$t_{j} = t_{6-j}$, shows that the new denominator is an {\em even} polynomial
in $w$. The symmetry of cosine about $\pi/2$ shows that the terms with odd 
power of $w$ have a vanishing integral. Thus, with $\psi = 2 \theta$, and
using the symmetry of the integrand to reduce the integral from $[0, \pi]$
to $[0, \pi/2]$, we obtain
\begin{equation}
U_{6} = \int_{0}^{\pi/2} \frac{r_{4} \cos^{4} \psi + r_{2} \cos^{2} \psi + 
r_{0} }{q_{6} \cos^{6} \psi + q_{4} \cos^{4} \psi + q_{2} \cos^{2} \psi + 
q_{0} } \, d \psi.  \nonumber 
\end{equation}
\noindent
The parameters $r_{j}, \, q_{j}$ have explicit formulas in terms of the 
original parameters of $U_{6}$. This even rational function of $\cos \psi$ 
can now be expressed in terms of $\cos  2 \psi$ to produce (letting 
$\theta \leftarrow 2 \psi$)
\begin{equation}
U_{6} = \int_{0}^{\pi} \frac{\alpha_{2} \cos^{2} \theta + \alpha_{1} 
\cos \theta +  \alpha_{0}}
{\beta_{3} \cos^{3} \theta + \beta_{2} \cos^{2} \theta + 
\beta_{1} \cos \theta +  \beta_{0} } \, d \theta. \nonumber 
\end{equation}
\noindent
The final change of variables $y = \tan \tfrac{\theta}{2}$ yields a new
rational form of the integrand:
\begin{equation}
U_{6} = \ift \frac{c_{1}y^{4} + d_{1}y^{2} + e_{1}}
{y^{6}+a_{1}y^{4}+b_{1}y^2+1} \, dy.
\label{usixi1}
\end{equation}
\noindent
Keeping track of the parameters, we have established: \\

\begin{theorem}
\label{landen6}
The integral 
\begin{equation}
U_{6} = \ift \frac{cx^{4} + dx^{2} + e}
{x^{6}+ax^{4}+bx^2+1} \, dx
\end{equation}
\noindent
is invariant under the change of parameters
\begin{eqnarray}
a_{1} & \leftarrow & \frac{ab + 5a+5b + 9}{(a+b+2)^{4/3}}, \label{ratlanden1} \\
b_{1} & \leftarrow & \frac{a + b + 6}{(a+b+2)^{2/3}}, \nonumber
\end{eqnarray}
\noindent
for the denominator parameters and 
\begin{eqnarray}
c_{1} & \leftarrow & \frac{c + d + e}{(a+b+2)^{2/3}}, \nonumber  \\
d_{1} & \leftarrow & \frac{(b+3)c + 2d + (a+3)e}{a+b+2}, \nonumber  \\
e_{1} & \leftarrow & \frac{c + e}{(a+b+2)^{1/3}} \nonumber
\end{eqnarray}
\noindent
for those of the numerator. 
\end{theorem}

\medskip

Theorem \ref{landen6} is the precise analogue of the elliptic  
Landen transformation (\ref{agm-1})
for the case of a rational integrand.  We call (\ref{ratlanden1}) a 
{\em rational Landen transformation}. This construction was 
first presented in \cite{boros2}.

\medskip

\subsection{Even rational Landen transformations}

More generally, there is a similar transformation
of coefficients for {\em any even rational integrand}; details appear in
\cite{boros1}.  We call these {\it even rational Landen Transformations}.  The obstruction in the general case comes from (\ref{polyinw});
one does not get a polynomial in $w = \cos 2 \theta$. \\

The method of proof for even rational integrals can be summarized 
as follows.  \\

\noindent
1) Start with an even rational integral:
\begin{equation}
U_{2p} = 
\ift \frac{\text{ even polynomial in } x }
{\text{ even polynomial in } x} \, dx.
\nonumber
\end{equation}

\medskip

\noindent
2) Symmetrize the denominator to produce 
\begin{equation}
U_{2p} = \ift \frac{ \text{ even polynomial in } x }{\text{ even  reciprocal
polynomial in } 
x } \, dx. 
\nonumber
\end{equation}
\noindent
The degree of the denominator is doubled.

\medskip

\noindent
3) Let $x = \tan \theta$. Then 
\begin{equation}
U_{2p} = \int_{0}^{\pi/2} \frac{\text{ polynomial in } \cos 2 \theta }
{\text{even polynomial in} \cos 2 \theta } \, d \theta. 
\nonumber
\end{equation}

\medskip

\noindent
4) Symmetry produced the vanishing of the integrands with an odd power of 
$\cos \theta$ in the numerator. We obtain
\begin{equation}
U_{2p} = \int_{0}^{\pi/2} \frac{\text{ even polynomial in } \cos 2 \theta  }
{\text{ even polynomial in } \cos 2 \theta  } \, d \theta. 
\nonumber
\end{equation}

\medskip

\noindent
5) Let $\psi = 2 \theta$ to produce 
\begin{equation}
U_{2p} = \int_{0}^{\pi} \frac{\text{ even polynomial in } \cos \psi }
{\text{ even polynomial in } \cos \psi } \, d \psi. 
\nonumber
\end{equation}
\noindent
Using symmetry this becomes an integral over $[0, \pi/2]$. 

\medskip

\noindent
6) Let $y = \tan \psi$ and use $\cos \psi = 1/\sqrt{1+y^{2}}$ to obtain
\begin{equation}
U_{2p} = \ift \frac{ \text{ even polynomial in } y }{\text{ even polynomial
in } y } \, dy. 
\nonumber
\end{equation}
\noindent
The degree of the denominator is half of what it was in Step 5. \\

Keeping track of the degrees one checks that the degree of the new rational 
function is the same as the original one, with new coefficients that appear
as functions of the old ones.

\section{A geometric interpretation} \label{sec-geometric} 
\setcounter{equation}{0}

We now present a geometric foundation of the general even rational 
Landen transformation (\ref{ratlanden1}) using the theory of 
Riemann surfaces. The text \cite{springer} provides an introduction to this 
theory, including definitions of objects we will refer to here.  The 
sequence of transformations
in section \ref{sec-ratlanden} can be achieved in one step by relating 
$\tan 2 \theta$ to $\tan \theta$. For historical reasons (this is what we 
did first) we present the details with {\em cotangent} instead of tangent. 

Consider the even rational integral 
\begin{equation}
I = \ift R(x) \, dx = \frac{1}{2} \int_{-\infty}^{\infty} R(x) \, dx. 
\label{int-1}
\end{equation}
\noindent
Introduce the new variable 
\begin{equation}
y = R_{2}(x) = \frac{x^{2}-1}{2x},
\end{equation}
\noindent
motivated by the identity $ \cot 2 \theta = R_{2}(\cot \theta)$. 
The function $R_{2}: \mathbb{R} \to \mathbb{R}$ is a two-to-one map. The 
sections of the inverse are
\begin{equation}
x = \sigma_{\pm}(y) = y \pm \sqrt{y^{2}+1}. 
\end{equation}
\noindent
Splitting the original integral as 
\begin{equation}
I = \int_{-\infty}^{0} R(x) \, dx + \int_{0}^{\infty} R(x) \, dx
\end{equation}
\noindent
and introducing $x = \sigma_{+}(y)$ in the first and $x = \sigma_{-}(y)$ 
in the second integral, yields
\begin{equation}
I = \int_{-\infty}^{\infty} \left( R_{+}(y) + R_{-}(y) \right) \, dy 
\label{int-2}
\end{equation}
\noindent
where 
\begin{equation}
R_{+}(y) = R( \sigma_{+}(y) ) + R(\sigma_{-}(y) )  \text{ and }
R_{-}(y) = \frac{y}{\sqrt{y^{2}+1} }
\left( R( \sigma_{+}(y) ) -R(\sigma_{-}(y) \right). 
\nonumber
\end{equation}
\noindent
A direct calculation shows that $R_{+}$ and $R_{-}$ are rational functions of 
degree at most that of $R$. 

The change of variables $y = R_{2}(x)$ converts the meromorphic differential
$\varphi = R(x) \, dx$ into 
\begin{eqnarray}
R( \sigma_{+}(y) ) \frac{d \sigma_{+}}{dy} + 
R( \sigma_{-}(y) ) \frac{d \sigma_{-}}{dy} & = & 
 \left( (R(\sigma_{+}) + R(\sigma_{-})) + 
\frac{y( R(\sigma_{+})-R(\sigma_{-}))}{\sqrt{y^{2}+1}}  \right) dy  \nonumber \\
& = & \left( R_{+}(y) + R_{-}(y) \right) \, dy. \nonumber 
\end{eqnarray}

The general situation is this: start with a finite ramified cover 
$\pi: X \to Y$ of Riemann surfaces and a 
meromorphic differential $\varphi$ 
on $X$. Let $U \subset Y$ be a 
simply connected domain that contains no critical 
values of $\pi$, and let $\sigma_{1}, \ldots, \sigma_{k}: U \to X$ be the 
distinct sections of $\pi$. Define
\begin{equation}
\pi_{*} \varphi\Big{|}_{U} = \sum_{j=1}^{k} \sigma_{j}^{*} \varphi. 
\end{equation}
\noindent
In \cite{hubbard1} we show that this construction preserves analytic
$1$-forms, that is, if $\varphi$ is an analytic $1$-form in $X$ then 
$\pi_{*} \varphi$ is an analytic $1$-form in $Y$. Furthermore, for any 
rectifiable curve $\gamma$ on $Y$, we have 
\begin{equation}
\int_{\gamma} \pi_{*} \varphi = \int_{\pi^{-1} \gamma} \varphi. 
\end{equation}
\noindent
In the case of projective space, this leads to the following: 

\begin{lemma}
\label{pi-rational}
If $\pi: \mathbb{P}^{1} \to \mathbb{P}^{1}$ is analytic, and $\varphi = 
R(z) \, dz$ with $R$ a rational function, then $\pi_{*} \varphi$ can be 
written as $R_{1}(z) \, dz$ with $R_{1}$ a rational function of degree at
most the degree of $R$. 
\end{lemma}

This is the generalization of the fact that the integrals in (\ref{int-1}) 
and (\ref{int-2}) are the same.

\section{A further generalization} \label{sec-generalization} 
\setcounter{equation}{0}

The procedure described in Section \ref{sec-ratlanden} can be extended with the
rational map $R_{m}$, defined by the identity
\begin{equation}
\cot m \theta  = R_{m}(\cot \theta).
\end{equation}
\noindent
Here $m \in \mathbb{N}$ is arbitrary greater or equal than $2$. We 
present some elementary 
properties of the rational function $R_{m}$. 

\begin{Prop}
The rational function $R_{m}$ satifies: \\

\noindent
1) For $m \in \mathbb{N}$ define 
\begin{equation}
P_{m}(x) := 
\sum_{j=0}^{\lfloor{ m/2 \rfloor}} (-1)^{j} \binom{m}{2j} x^{m-2j}
\text{ and } 
Q_{m}(x)  :=  \sum_{j=0}^{\lfloor{ (m-1)/2 \rfloor}} 
(-1)^{j} \binom{m}{2j+1} x^{m-(2j+1)}. \nonumber 
\end{equation}
\noindent
Then $R_{m} := P_{m}/Q_{m}$. \\

\noindent
2) The function $R_{m}$ is 
conjugate to $ \, f_{m}(x) := x^{m}$ via $ \, M(x) := 
\frac{x+i}{x-i}$; that is, $R_{m} = M^{-1} \circ f_{m} \circ M$. \\

\noindent
3) The polynomials $P_{m}$ and $Q_{m}$ have simple real zeros given by
\begin{equation}
p_{k} := \cot \left( \frac{(2k+1) \pi}{2m} \right) \text{ for } 
0 \leq k \leq m-1, \text{ and }
q_{k} := \cot \left( \frac{k \pi}{m} \right) \text{ for } 
1 \leq k \leq m-1. \nonumber 
\end{equation}
\end{Prop}

If we change the domain to the entire real line, we can, using the rational 
substitutions $R_m(x) \mapsto x$, produce a rational Landen transformation for 
an arbitrary integrable rational function $R(x) = B(x)/A(x)$ for each 
integer value of $m$. The 
result is a new 
list of coefficients, from which one produces a second rational function
$R^{(1)}(x) = J(x)/H(x)$ with
\begin{equation}
\int_{-\infty}^{\infty} \frac{B(x)}{A(x)} \, dx = 
\int_{-\infty}^{\infty} \frac{J(x)}{H(x)} \, dx.
\end{equation}
\noindent
Iteration of this procedure yields a sequence 
$\mathbf{x}_{n}$, that has a limit $\mathbf{x}_{\infty}$ with convergence 
of order $m$, that is, 
\begin{equation}
\| \mathbf{x}_{n+1} - \mathbf{x}_{\infty} \| \leq C 
\| \mathbf{x}_{n} - \mathbf{x}_{\infty} \|^{m}.
\end{equation}
\noindent
We describe this procedure here in the 
form of an algorithm; proofs appear in \cite{manna-moll2}.  \\

Lemma \ref{pi-rational} applied to the map 
$\pi(x) = R_{m}(x)$, viewed as ramified cover of $\mathbb{P}^{1}$, 
guarantees the existence of a such new rational function $R^{(1)}$.
The question of effective computation of the coefficients
of $J$ and $H$ is discussed below. In particular, we show that all these 
calculations can be done symbolically. \\

\noindent $\bullet$ {\bf Algorithm for 
Deriving Rational Landen Transformations}

\medskip

\noindent
{\bf Step 1}. The initial data is a rational function $R(x) := B(x)/A(x)$. We 
assume that $A$ and $B$ are polynomials with real coefficients and $A$ has 
no real zeros and write 
\begin{equation}
A(x) := \sum_{k=0}^{p} a_{k}x^{p-k} \text{ and } 
B(x) := \sum_{k=0}^{p-2} b_{k}x^{p-2-k}.
\end{equation}

\medskip

\noindent
{\bf Step 2}. Choose a positive integer $m \geq 2$. \\


\noindent
{\bf Step 3}. Introduce the polynomial
\begin{equation}
H(x) := \text{Res}_{z}( A(z), P_{m}(z) - x Q_{m}(z))
\end{equation}
\noindent
and write it as 
\begin{equation}
H(x) := \sum_{l=0}^{p} e_{l}x^{p-l}.
\end{equation}
\noindent
The polynomial $H$ is thus defined as the determinant of the Sylvester matrix 
which is formed of the polynomial coefficients. As such, the coefficients 
$e_{l}$ of $H(x)$ themselves are integer polynomials in the $a_{i}$. Explicitly,
\begin{equation}
e_{l} = (-1)^{l} a_{0}^{m} \prod_{j=1}^{p} Q_{m}(x_{j}) \times 
\sigma_{l}^{(p)}(R_{m}(x_{1}), \, R_{m}(x_{2}), \, \ldots, R_{m}(x_{p}) ),
\end{equation}
\noindent
where $\{ x_{1}, x_{2}, \cdots, x_{p} \}$ are the roots of $A$, each
written according to multiplicity. The functions 
$\sigma_{l}^{(p)}$ are the elementary symmetric functions in $p$ variables 
defined by
\begin{equation}
\prod_{l=1}^{p} (y - y_{l}) = \sum_{l=0}^{p} (-1)^{l} 
\sigma_{l}^{(p)}(y_{1}, \cdots, y_{p}) y^{p-l}. 
\end{equation}
\noindent
It is possible to compute the coefficients $e_{l}$ symbolically
from the coefficients of $A$, without the knowledge of the roots 
of $A$. 

Also define 
\begin{equation}
E(x) := H(R_{m}(x)) \times Q_{m}(x)^{p}.
\end{equation}

\medskip

\noindent
{\bf Step 4}. The polynomial $A$ divides $E$ and we denote the quotient by 
$Z$. The coefficients of $Z$ are integer polynomials in the $a_{i}$. \\

\noindent
{\bf Step 5}.  Define the polynomial $C(x) := B(x) Z(x)$. \\

\noindent
{\bf Step 6}. There exists a polynomial $J(x)$, whose coefficients have an 
explicit formula in terms of the coefficients $c_{j}$ of $C(x)$, such
that
\begin{equation}
\int_{-\infty}^{\infty} \frac{B(x)}{A(x)} \,dx  = 
\int_{-\infty}^{\infty} \frac{J(x)}{H(x)} \,dx. 
\end{equation}
\noindent
This new integrand is the rational function whose existence is guaranteed by
Lemma \ref{pi-rational}. The
explicit computation of the coefficients of $J$ can be found 
in \cite{manna-moll2}. This 
is the {\em rational Landen transformation} of order $m$. 

\medskip

\noindent
\begin{example}\label{simpleex}Completing the algorithm with $m=3$ and the 
rational function
\begin{equation}
R(x) = \frac{1}{ax^{2} +bx +c },
\end{equation}
\noindent
produces the result stated below. Notice that the values of the iterates 
are ratios of integer polynomials of degree 3, as was stated above.  The 
details of this example appear in \cite{manna-moll1}.

\begin{Thm}
\label{thm-monthly}
The integral
\begin{equation}
I = \int_{-\infty}^{\infty} \frac{dx}{ax^{2} + bx + c} 
\end{equation}
\noindent
is invariant under the transformation 
\begin{eqnarray}
a & \mapsto & \frac{a}{\Delta} \left( (a+3c)^{2} - 3b^{2} \right),
 \label{quadratic3}  \\
b & \mapsto & \frac{b}{\Delta} \left( 3(a-c)^{2} - b^{2} \right), \nonumber \\
c & \mapsto & \frac{c}{\Delta} \left( (3a+c)^{2} - 3b^{2} \right), \nonumber
\end{eqnarray}
\noindent
where $\Delta := (3a+c)(a+3c)-b^{2}$. 
The condition $b^{2}-4ac < 0$, imposed to ensure convergence of the integral,
is preserved by the iteration.
\end{Thm}
\end{example}

\medskip

\begin{example}
In this example we follow the steps described above in order to produce a
rational Landen transfromation of order $2$ for the integral

\begin{equation}
I = \ifft \frac{b_{0}x^{4} + b_{1}x^{3} + b_{2}x^{2} + b_{3}x + b_{4}}
{a_{0}x^{6} + a_{1}x^{5} + a_{2}x^{4} + a_{3}x^{3} + a_{4}x^{2} + a_{5}x 
+ a_{6}} \, dx. 
\end{equation}

\medskip

Recall that the algorithm starts with a rational function $R(x)$ and 
produces a new function ${\mathfrak{L}}_{2}(R(x))$ satisfying 
\begin{equation}
\ifft R(x) \, dx = \ifft {\mathfrak{L}}_{2}(R(x)) \, dx. \label{equalint}
\end{equation}

\medskip

\noindent
{\bf Step 1}. The initial data is $R(x) = B(x)/A(x)$ with 
\begin{equation}
A(x) = a_{0}x^{6} + a_{1}x^{5} + a_{2}x^{4} + a_{3}x^{3} + a_{4}x^{2} + a_{5}x 
+ a_{6}, 
\end{equation}
\noindent
and 
\begin{equation}
B(x) = b_{0}x^{4} + b_{1}x^{3} + b_{2}x^{2} + b_{3}x + b_{4}. 
\end{equation}
\noindent
The parameter $p$ is the degree of $A$, so $p=6$. \\

\noindent
{\bf Step 2}. We choose $m=2$ to  produce a method of order $2$. The 
algorithm employs the polynomials $P_{2}(z) = z^{2} -1$ and $Q_{2}(z) = 2z$. \\

\noindent
{\bf Step 3}. The polynomial 
\begin{equation}
H(x) := \text{Res}_{z}( A(z), z^{2}-1 - 2xz)
\end{equation}
\noindent
is computed with the Mathematica command \texttt{Resultant} to obtain
\begin{equation}
H(x) = e_{0}x^{6} + e_{1}x^{5} + e_{2}x^{4} + e_{3}x^{3} + e_{4}x^{2} 
+ e_{5}x + e_{6}, 
\label{newdenom}
\end{equation}
\noindent
with
\begin{eqnarray}
 & & \label{valuesofe} \\
e_{0} & = & 64a_{0}a_{6}, \nonumber  \\
e_{1} & = & -32(a_{0}a_{5}-a_{1}a_{6}), \nonumber \\
e_{2} & = & 16(a_{0}a_{4}-a_{1}a_{5}+6a_{0}a_{6}+a_{2}a_{6}), \nonumber \\
e_{3} & = & -8(a_{0}a_{3}-a_{1}a_{4}+5a_{0}a_{5}+a_{2}a_{5}-5a_{1}a_{6}-a_{3}a_{6}), \nonumber \\
e_{4}& = & 4(a_{0}a_{2}-a_{1}a_{3}+4a_{0}a_{4}+a_{2}a_{4}-4a_{1}a_{5}-a_{3}a_{5}+9a_{0}a_{6}+4a_{2}a_{6}+a_{4}a_{6}), \nonumber \\
e_{5} & = & -2(a_{0}a_{1}-a_{1}a_{2}+3a_{0}a_{3}+a_{2}a_{3}-3a_{1}a_{4}-a_{3}a_{4}+5a_{0}a_{5}+ \nonumber \\
 & & + 3a_{2}a_{5}+a_{4}a_{5}-5a_{1}a_{6}-3a_{3}a_{6}-a_{5}a_{6}), 
\nonumber \\
e_{6} & = & (a_{0}-a_{1}+a_{2}-a_{3}+a_{4}-a_{5}+a_{6})(a_{0}+a_{1}+a_{2}+a_{3}+a_{4}+a_{5}+a_{6}). \nonumber 
\end{eqnarray}

\noindent
The polynomial $H(x)$ is the denominator of the integrand 
${\mathfrak{L}}_{2}(R(x))$ in 
(\ref{equalint}). \\

In Step 3 we also define
\begin{equation}
E(x) = H(R_{2}(x)) Q_{2}^{6}(x) = H \left( \frac{x^{2}-1}{2x} \right) 
\cdot (2x)^{6}. 
\label{Hdef}
\end{equation}
\noindent
The function $E(x)$ is a polynomial of degree $12$, written as 
\begin{equation}
E(x) = \sum_{k=0}^{12} \alpha_{k}x^{12-k}.
\end{equation}
\noindent
Using the expressions for $e_{j}$ in (\ref{valuesofe}) in (\ref{Hdef}) yields
\begin{eqnarray}
\alpha_{0} & = & \alpha_{12} = 64a_{0}a_{6} \label{valuesofalpha} \\
\alpha_{1} & = & - \alpha_{11} = -64(a_{0}a_{5}-a_{1}a_{6}), \nonumber \\
\alpha_{2} & = & \alpha_{10} = 64(a_{0}a_{4}-a_{1}a_{5}+a_{2}a_{6}), 
\nonumber \\
\alpha_{3} & = & - \alpha_{9} = -64(a_{0}a_{3}-a_{1}a_{4}+a_{2}a_{5}-a_{3}a_{6}), \nonumber \\
\alpha_{4} & =  & \alpha_{8} = 64(a_{0}a_{2}-a_{1}a_{3}+a_{2}a_{4}-a_{3}a_{5}+a_{4}a_{6}), \nonumber \\
\alpha_{5} & = & - \alpha_{7} = -64(a_{0}a_{1}-a_{1}a_{2}+a_{2}a_{3}-a_{3}a_{4}+a_{4}a_{5}-a_{5}a_{6}), \nonumber \\
\alpha_{6} & = & 64(a_{0}^{2}-a_{1}^{2}+a_{2}^{2}-a_{3}^{2}+a_{4}^{2}-a_{5}^{2}+a_{6}^{2}). \nonumber 
\end{eqnarray}

\medskip

\noindent
{\bf Step 4}. The polynomial $A(x)$ always divides $E(x)$. The quotient is
denoted by $Z(x)$. The values of $\alpha_{j}$ given in (\ref{valuesofalpha}) 
produce
\begin{equation}
Z(x) = 64(a_{0}-a_{1}x+a_{2}x^{2}-a_{3}x^{3}+a_{4}x^{4}-a_{5}x^{5}+a_{6}x^{6}).
\end{equation} 

\medskip

\noindent
{\bf Step 5}. Define the polynomial $C(x) := B(x)Z(x)$. In this case, $C$ is
of degree $10$, written as 
\begin{equation}
C(x) = \sum_{k=0}^{10} c_{k}x^{10-k}, 
\end{equation}
\noindent
and the coefficients $c_{k}$ are given by
\begin{eqnarray}
& & \label{coefficientsc} \\
c_{0} & = & 64a_{6}b_{0}, \nonumber \\
c_{1} & = & -64(a_{5}b_{0}-a_{6}b_{1}), \nonumber \\
c_{2} & = & 64(a_{4}b_{0}-a_{5}b_{1}+a_{6}b_{2}), \nonumber \\
c_{3} & = & -64(a_{3}b_{0}-a_{4}b_{1}+a_{5}b_{2}-a_{6}b_{3}), \nonumber \\
c_{4} & = & 64(a_{2}b_{0}-a_{3}b_{1}+a_{4}b_{2}-a_{5}b_{3}+a_{6}b_{4}), \nonumber \\
c_{5} & = & -64(a_{1}b_{0}-a_{2}b_{1}+a_{3}b_{2}-a_{4}b_{3}+a_{5}b_{4}), 
\nonumber \\
c_{6} & = & 64(a_{0}b_{0}-a_{1}b_{1}+a_{2}b_{2}-a_{3}b_{3}+a_{4}b_{4}),
\nonumber \\
c_{7} & = & 64(a_{0}b_{1}-a_{1}b_{2}+a_{2}b_{3}-a_{3}b_{4}),
\nonumber \\
c_{8} & = & 64(a_{0}b_{2}-a_{1}b_{3}+a_{2}b_{4}),
\nonumber \\
c_{9} & = & 64(a_{0}b_{3}-a_{1}b_{4}), 
\nonumber \\
c_{10} & = & 64a_{0}b_{4}. \nonumber 
\end{eqnarray}

\medskip

\noindent
{\bf Step 6} produces the numerator $J(x)$ of the new integrand
${\mathfrak{L}_{2}}(R(x))$ from  the coefficients $c_{j}$ given in 
(\ref{coefficientsc}). The function $J(x)$ is a polynomial of degree $4$,
written as
\begin{equation}
J(x) = \sum_{k=0}^{4}j_{k}x^{4-k}.
\label{newnumer}
\end{equation}
\noindent
Using the values of (\ref{coefficientsc}) we obtain
\begin{eqnarray}
& & \label{coefficientsj} \\
j_{0} & = & 32(a_{6}b_{0}+a_{0}b_{4}), \nonumber \\
j_{1} & = & -16(a_{5}b_{0} -a_{6}b_{1}+a_{0}b_{3} - a_{1}b_{4}), 
\nonumber \\
j_{2} & = & 8(a_{4}b_{0}+3a_{6}b_{0}-a_{5}b_{1}+a_{0}b_{2}+a_{6}b_{2}-a_{1}b_{3}+3a_{0}b_{4}+a_{2}b_{4}),
\nonumber \\
j_{3} & = & -4(a_{3}b_{0}+2a_{5}b_{0}+a_{0}b_{1}-a_{4}b_{1}-2a_{6}b_{1}-a_{1}b_{2}+a_{5}b_{2}+ \nonumber \\
& & + 2a_{0}b_{3}+a_{3}b_{3}-a_{6}b_{3}-2a_{1}b_{4}-a_{3}b_{4}),
\nonumber \\
j_{4} & = & 2( a_{0}b_{0}+a_{2}b_{0}+a_{4}b_{0}+a_{6}b_{0}-a_{1}b_{1}-a_{3}b_{1}-a_{5}b_{1}+a_{0}b_{2}+ \nonumber \\
& & + a_{2}b_{2}+a_{4}b_{2}+a_{6}b_{2}-a_{1}b_{3}-a_{3}b_{3}-a_{5}b_{3}+a_{0}b_{4}+a_{2}b_{4}+a_{4}b_{4}+a_{6}b_{4}). 
\nonumber 
\end{eqnarray}
\noindent
The explicit computation of the coefficients of $J$ can be found 
in \cite{manna-moll2}. 

\medskip

The new rational function is 
\begin{equation}
{\mathfrak{L}}_{2}(R(x)) := \frac{J(x)}{H(x)},
\end{equation}
\noindent
with $J(x)$ given in (\ref{newnumer}) and $H(x)$ in (\ref{newdenom}). The 
transformation is 
\begin{equation}
\frac{b_{0}x^{4} + b_{1}x^{3} + b_{2}x^{2} + b_{3}x + b_{4}}
{a_{0}x^{6} + a_{1}x^{5} + a_{2}x^{4} + a_{3}x^{3} + a_{4}x^{2} + a_{5}x 
+ a_{6}} \mapsto 
\frac{j_{0}x^{4} + j_{1}x^{3} + j_{2}x^{2} + j_{3}x + j_{4}}
{e_{0}x^{6} + e_{1}x^{5} + e_{2}x^{4} + e_{3}x^{3} + e_{4}x^{2} + e_{5}x 
+ e_{6}}. \nonumber
\end{equation}
\noindent
The numerator coefficients are given in (\ref{valuesofe}) and the denominator
ones in (\ref{coefficientsj}), explicitly as polynomials in the 
coefficients of the original rational function. The generation of these
polynomials is a completely symbolic procedure. 

\medskip

The first two steps of this algorithm, applied to the definite integral
\begin{equation}
\ifft \frac{dx}{x^{6}+x^{3}+1} = 
\frac{\pi}{9} \left( 2 \sqrt{3} \cos(\pi/9) + \sqrt{3} \cos(2 \pi/9) + 
3 \sin(2 \pi/9) \right), 
\end{equation}
\noindent
produces the identities
\begin{eqnarray}
\ifft \frac{dx}{x^{6}+x^{3}+1}  & = & 
\ifft \frac{2(16x^4+12x^2+2x+2)}{64x^6 + 96x^4+36x^2+3} \, dx \nonumber \\
& = & 
\ifft \frac{4(2816x^4-1024x^3+8400x^2-884x+5970)}
{12288x^6+59904x^4+87216x^2+39601} \, dx. \nonumber
\end{eqnarray}

\end{example}

\medskip

The convergence of the iterations of rational Landen transformations is 
discussed in the next section.

\section{The issue of  convergence} \label{sec-convergence} 
\setcounter{equation}{0}

The convergence of the double sequence 
$(a_{n},b_{n})$ appearing in the elliptic Landen transformation (\ref{agm-1})
is easily established. Assume $0 < b_{0} \leq a_{0}$, then the
arithmetic-geometric inequality yields
$b_{n} \leq b_{n+1} \leq a_{n+1} \leq a_{n}$. Also 
\begin{equation}
0 \leq a_{n+1}-b_{n+1} = \frac{1}{2} \frac{(a_{n}-b_{n})^{2}}
{ (\sqrt{a_{n}} + \sqrt{b_{n}})^{2}}. 
\label{conv-elliptic}
\end{equation}
\noindent
This shows that $a_{n}$ and $b_{n}$ have a common limit: $M = \agm(a,b)$, the 
arithmetic-geometric of $a$ and $b$. The convergence is quadratic:
\begin{equation}
| a_{n+1} - M | \leq C | a_{n} - M |^{2}, 
\end{equation}
\noindent
for some constant  $C>0$ independent of $n$.  Details can be 
found in \cite{borwein1}.

The Landen transformations produce maps on the space 
of coefficients of the integrand.  
In this section, we discuss the convergence of the rational Landen 
transformations. This discussion is divided in two cases: \\

\noindent
{\bf Case 1: the half-line}. Let $R(x)$ be an even rational function, written
as $R(x) = P(x)/Q(x)$, with 
\begin{equation}
P(x) =  \sum_{k=0}^{p-1} b_{k}x^{2(p-1-k)}, \, 
Q(x) =  \sum_{k=0}^{p} a_{k}x^{2(p-k)},
\end{equation}
\noindent
and $a_{0}=a_{p}=1$. The {\em parameter space} is 
\begin{equation}
\mathfrak{P}_{2p}^{+} = \{ (a_{1}, \cdots, a_{p-1}; b_{0}, \cdots, b_{p-1}) \}
\subset \mathbb{R}^{p-1} \times \mathbb{R}^{p}. 
\end{equation}
\noindent
We write 
\begin{equation}
\mathbf{a}:= (a_{1}, \cdots, a_{p-1}) \text{ and } 
\mathbf{b}:= (b_{0}, \cdots, b_{p}). 
\end{equation}

Define 
\begin{equation}
\Lambda_{2p} = \{ (a_{1}, \cdots, a_{p-1}) \in \mathbb{R}^{p-1}: 
\ift R(x) \, dx \text{ is finite} \, \}. 
\label{lambda}
\end{equation}
\noindent
Observe that the convergence of the integral depends only on the parameters in
the denominator. 

The Landen transformations provide a map
\begin{equation}
\Phi_{2p}: \mathfrak{P}_{2p}^{+} \to \mathfrak{P}_{2p}^{+}
\end{equation}
\noindent
that preserves the integral. Introduce the notation
\begin{equation}
\mathbf{a}_{n} = (a_{1}^{(n)}, \cdots, a_{p-1}^{(n)}) \text{ and }
\mathbf{b}_{n} = (b_{0}^{(n)}, \cdots, b_{p}^{(n)}),
\end{equation}
where 
\begin{equation}
(\mathbf{a}_{n}, \mathbf{b}_{n}) =  
\Phi_{2p}(\mathbf{a}_{n-1}, \mathbf{b}_{n-1})
\end{equation}
\noindent
are the iterates of the map $\Phi_{2p}$. \\

The result that one expects is this: 

\begin{Thm}
\label{converg-1}
The region $\Lambda_{2p}$ is invariant under the map $\Phi_{2p}$.  Moreover
\begin{equation}
\mathbf{a}_{n} \to \Bigg( {\tiny \binom{p}{1}, \binom{p}{2}, \cdots, \binom{p}{p-1}}
\Bigg), 
\end{equation}
\noindent
and there exists a number $L$, that depends on the initial conditions, such
that
\begin{equation}
\mathbf{b}_{n} \to \Bigg( {\tiny \binom{p-1}{0}L, \binom{p-1}{1}L, \cdots, 
\binom{p-1}{p-1} L} \Bigg).
\end{equation}
\noindent
This is equivalent to say that the sequence of rational functions formed by
the Landen transformations, converge to $L/(x^{2}+1)$. 
\end{Thm}

This was established in \cite{hubbard1} using the geometric language of
Landen transformations which, while unexpected, is satisfactory.

\begin{Thm}\label{converg-2}
Let $\varphi$ be a $1$-form, holomorphic in a neighborhood of 
$\mathbb{R} \subset \mathbb{P}^{1}$. Then 
\begin{equation}
\lim\limits_{n \to \infty} (\pi_{*})^{n} \varphi = 
\frac{1}{\pi} \left( \int_{-\infty}^{\infty} \varphi \right) \frac{dz}{1+z^{2}},
\end{equation}
\noindent
where the convergence is uniform on compact subsets of $U$, the neighborhood 
in the definition of $\pi_{*}$.
\end{Thm}
\noindent
The proof is detailed for the map $\pi(z) = \frac{z^{2}-1}{2z} = R_{2}(z)$, 
but it extends without difficulty to the generalization $R_{m}$. \\ 

Theorem \ref{converg-2} can be equivalently reformulated as:

\begin{Thm}
The iterates of the Landen transformation starting at 
$(\mathbf{a}_{0}, \mathbf{b}_{0}) \in \mathfrak{P}_{2p}^{+}$ 
converge (to the limit stated in 
Theorem \ref{converg-1} ) if and only if the integral formed by the 
initial data is finite.
\end{Thm}

It would be desirable to establish this result by purely dynamical techniques.
This has been established only for the case $p=3$. In that case the Landen
transformation for
\begin{equation}
U_{6} := \ift \frac{cx^{4} + dx^{2} + e}{x^{6}+ax^{4}+bx^{2}+1} \, dx
\end{equation}
\noindent
is 
\begin{eqnarray}
a_{1} 
& \leftarrow & \frac{ab + 5a+5b + 9}{(a+b+2)^{4/3}}, \label{ratlanden1a} \\
b_{1} & \leftarrow & \frac{a + b + 6}{(a+b+2)^{2/3}}, \nonumber
\end{eqnarray}
\noindent
coupled with 
\begin{eqnarray}
c_{1} & \leftarrow & \frac{c + d + e}{(a+b+2)^{2/3}}, \nonumber  \\
d_{1} & \leftarrow & \frac{(b+3)c + 2d + (a+3)e}{a+b+2}, \nonumber  \\
e_{1} & \leftarrow & \frac{c + e}{(a+b+2)^{1/3}}. \nonumber
\end{eqnarray}
\noindent
The region 
\begin{equation}
\Lambda_{6} = \{ (a,b) \in \mathbb{R}^{2}:  U_{6} < \infty \}
\end{equation}
\noindent
is described by the discriminant curve $\mathfrak{R}$, the zero set of the
polynomial 
\begin{equation}
R(a,b) = 4a^{3}+4b^{3} - 18ab -a^{2}b^{2} + 27.
\end{equation}

\input{epsf}

\begin{figure}
\epsfysize=3 in
\epsfxsize=3 in
\epsffile{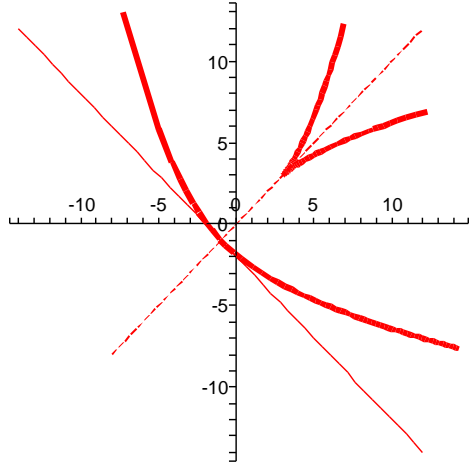}
\caption{Zero locus of $R(a,b)$}
\label{locus}
\end{figure}

\noindent
This zero set, shown in Figure \ref{locus}, has 
two connected components: the first 
one $\mathfrak{R}_{+}$ 
contains $(3,3)$ as a cusp and the second one $\mathfrak{R}_{-}$, given by 
$R_{-}(a,b) = 0$, is disjoint from the first quadrant. The branch 
$\mathfrak{R}_{-}$ is the boundary of the set $\Lambda_{6}$. 

The identity 
\begin{equation}
R(a_{1},b_{1}) = \frac{(a-b)^{2} \, R(a,b) }{(a+b+2)^{4}},
\end{equation}
shows that $\partial \mathfrak{R}$ is invariant under $\Phi_{6}$. By examining 
the effect of this map along lines of slope $-1$, we obtain a direct 
parametrization of the flow {\em on} the discriminant curve. Indeed, this 
curve is parametrized by 
\begin{equation}
a(s) = \frac{s^{3}+4}{s^{2}} \text{ and } b(s) = \frac{s^{3}+16}{4s}. 
\end{equation}
\noindent
Then 
\begin{equation}
\varphi(s) = \left( \frac{4(s^{2}+4)^{2}}{s(s+2)^{2}} \right)^{1/3}
\end{equation}
\noindent
gives the image of the Landen transformation $\Phi_{6}$; that is,
\begin{equation}
\Phi_{6}(a(s),b(s)) = ( a (\varphi(s)), b(\varphi(s) ). 
\end{equation}

The map $\Phi_{6}$ has three fixed points: $(3,3),$ that is super-attracting, 
a saddle point $P_{2}$ on the lower branch 
$\mathfrak{R}_{-}$ of the discriminant curve, and
a third unstable spiral below
this lower branch. In \cite{marc-moll} we prove:

\begin{Thm}
The lower branch of the discriminant curve is the curve $\Lambda_{6}$. 
This curve is also the global unstable manifold of the saddle point 
$P_{2}$. Therefore the iterations of $\Phi_{6}$ starting at $(a,b)$ 
converge if and only if 
the integral $U_{6}$, formed with the parameters $(a,b)$, is 
finite. Moreover, 
$(a_{n},b_{n}) \to (3,3)$ quadratically and there exists a number $L$ such that
$(c_{n}, d_{n}, e_{n}) \to (1,2,1)L$. 
\end{Thm}

\medskip

The next result provides an analogue of the $\agm$ (\ref{agmlaw}) for the 
rational case.  The main differences here are that our iterates converge to 
an {\it algebraic} number and we achieve {\it order-m} convergence. \\

\noindent
{\bf Case 2: The whole-line}: This works for any choice of positive integer $m$. Let 
$R(x)$ be a rational function, written as 
$R(x) = B(x)/A(x)$. Assume that the coefficients of $A$ and $B$ are real, that
$A$ has no real zeros and that $\text{deg}(B) \leq \text{deg}(A)-2$. These 
conditions are imposed to guarantee the existence of 
\begin{equation}
I = \int_{-\infty}^{\infty} R(x) \, dx. 
\end{equation}
\noindent
In particular $A$ must be of even degree, and we write
\begin{equation}
A(x) = \sum_{k=0}^{p} a_{k}x^{p-k} \text{ and } 
B(x) = \sum_{k=0}^{p-2} b_{k}x^{p-2-k}. 
\end{equation}
\noindent
We can also require that $deg(gcd(A,B)) = 0$. 

The class of such rational functions will be denoted by ${\mathfrak{R}}_{p}$. 
\\

The algorithm presented in Section \ref{sec-generalization}
provides a transformation on  the parameters
\ba
{\mathfrak{P}}_{p} & := & \{ a_{0}, \, a_{1}, \cdots, a_{p}; \, b_{0}, \,
b_{1}, \cdots, b_{p-2} \, \} = \mathbb{R}^{p+1} \times \mathbb{R}^{p-1}
\label{initialpara}
\ea
\no
of $R \in {\mathfrak{R}}_{p}$ that preserves the integral $I$. In fact, we
produce a family of maps, indexed by $m \in \mathbb{N}$,
$$
{\mathfrak{L}}_{m,p}: {\mathfrak{R}}_{p} \to {\mathfrak{R}}_{p},
$$
\no
such that
\ba
\ifft R(x) \, dx & = & \ifft {\mathfrak{L}}_{m,p}(R(x)) \, dx.
\ea
\no
The maps ${\mathfrak{L}}_{m,p}$ induce a {\em rational Landen
transformation} on the parameter space:
\ba
\Phi_{m,p}: {\mathfrak{P}}_{p} \to {\mathfrak{P}}_{p}
\ea
\noindent
by simply listing the coefficients of the 
function ${\mathfrak{L}}_{m,p}(R(x))$.

The original integral is written in the form 
\ba I & = &
\frac{b_{0}}{a_{0}} \ifft \frac{x^{p-2} + b_{0}^{-1}b_{1}x^{p-3} +
b_{0}^{-1}b_{2}x^{p-4} + \cdots + b_{0}^{-1}b_{p-2}} {x^{p} +
a_{0}^{-1}a_{1}x^{p-1} + a_{0}^{-1}a_{2}x^{p-2} + \cdots +
a_{0}^{-1}a_{p}} \, dx. 
\ea 
\no 
The Landen transformation
generates a sequence of coefficients, 
\ba {\mathfrak{P}}_{p,n} & :=
& \{ a_{0}^{(n)}, \, a_{1}^{(n)}, \cdots, a_{p}^{(n)}; \,
b_{0}^{(n)}, \, b_{1}^{(n)}, \cdots, b_{p-2}^{(n)} \, \}~, \ea \no
with ${\mathfrak{P}}_{p,0} = {\mathfrak{P}}_{p}$ as in
(\ref{initialpara}). We expect that, as $n \to \infty$,
\begin{equation}
{\mathbf{x}}_{n} := \left( \frac{a_{1}^{(n)}}{a_{0}^{(n)}},
\frac{a_{2}^{(n)}}{a_{0}^{(n)}},  \cdots,
\frac{a_{p}^{(n)}}{a_{0}^{(n)}},
\frac{b_{1}^{(n)}}{b_{0}^{(n)}},
\frac{b_{2}^{(n)}}{b_{0}^{(n)}},  \cdots,
\frac{b_{p-2}^{(n)}}{b_{0}^{(n)}} \right)
\end{equation}
\no
converges to
\begin{equation}
{\mathbf{x}}_{\infty} := \Bigg( 0, \binom{q}{1},0, \binom{q}{2}, \cdots,
\binom{q}{q}; 0, \binom{q-1}{1},0, \binom{q-1}{2}, \cdots,
\binom{q-1}{q-1} \Bigg)~,
\end{equation}
\no
where $q = p/2$. Moreover, we should have
\begin{equation}
\| {\mathbf{x}}_{n+1} - {\mathbf{x}}_{\infty} \| \leq C 
\| {\mathbf{x}}_{n} - {\mathbf{x}}_{\infty} \|^{m}. 
\label{order-m}
\end{equation}
\noindent
The invariance of the integral then shows that
\ba
\frac{b_{0}^{(n)}}{a_{0}^{(n)}} \to \frac{1}{\pi} I.
\label{conv-int}
\ea
\noindent

This produces an iterative method to evaluate the integral of a rational 
function. The method's convergence is of order-$m$. \\

The convergence of these iterations, and in particular the 
bound (\ref{order-m}), can 
be established by the argument presented in 
section \ref{sec-geometric}. Thus, the transformation
${\mathfrak{L}}_{m,p}$ leads to a sequence that has order-$m$ convergence. 
We expect to develop these ideas into an efficient numerical method for
integration.  \\

\medskip

We choose to measure the convergence of a sequence of vectors to $0$ is in 
the $L_{2}-$norm,
\ba \Vert v \Vert_{2} & = & \frac{1}{\sqrt{2p-2}}
\left( \sum_{k=1}^{2p-2} \| v_{k} \|^{2} \right)^{1/2}, \ea \no and
also the $L_{\infty}$-norm, \ba \Vert v \Vert_{\infty} & = &
\text{Max} \left\{ \| v_{k} \|: 1 \leq k \leq 2p-2 \, \right\}.
\ea \no The rational functions appearing as integrands have
rational coefficients, so, as a measure of their complexity, we
take the largest number of digits of these coefficients. This
appears in the column marked {\em size}.

\bigskip

The following tables illustrate the iterates of rational Landen
transformations of order $2, \, 3$ and $4$, applied to the example
\ba 
I & = & \ifft \frac{3x+5}{x^{4}+ 14x^3+74x^2+184x+208} \, dx =
- \frac{7 \pi}{12}. \nonumber 
\ea

\medskip

\no The first column gives the $L_{2}$-norm of $u_{n} -
u_{\infty}$, the second its $L_{\infty}$-norm, the third presents
the relative error in (\ref{conv-int}), and in the last column we
give the size of the rational integrand. At each step, we verify
that the new rational function integrates to $-7 \pi/12$.

\begin{center}
Method of order $2$
\end{center}

\medskip

\begin{center}
\begin{tabular}{||c|c|c|c|c||}
\hline
$n$ & $L_{2}$-norm & $L_{\infty}$-norm & Error & Size \\ \hline
$1$ & $58.7171$ & $69.1000$ & $1.02060$ & $5$ \\
$2$ & $7.444927$ & $9.64324$ & $1.04473$ & $10$ \\
$3$ & $4.04691$ & $5.36256$ & $0.945481$ & $18$ \\
$4$ & $1.81592$ & $2.41858$ & $1.15092$ & $41$ \\
$5$ & $0.360422$ & $0.411437$ & $0.262511$ & $82$ \\
$6$ & $0.0298892$ & $0.0249128$ & $0.0189903$ & $164$ \\
$7$ & $0.000256824$ & $0.000299728$ & $0.0000362352$ & $327$ \\
$8$ & $1.92454 \times 10^{-8}$ & $2.24568 \times 10^{-8}$ & $1.47053 \times 10^{-8}$ & $659$ \\
$9$ & $1.0823 \times 10^{-16}$ & $1.2609 \times 10^{-16}$ & $8.2207 \times 10^{-17}$
& $1318$ \\
\hline
\end{tabular}
\end{center}

\medskip

As expected, for the method of order $2$,  we 
observe quadratic convergence in the $L_{2}-$norm
and also in the $L_{\infty}-$norm. The size of the coefficients of the 
integrand is approximately doubled at each iteration. \\

\begin{center}
Method of order $3$
\end{center}

\medskip

\begin{center}
\begin{tabular}{||c|c|c|c|c||}
\hline
$n$ & $L_{2}$-norm & $L_{\infty}$-norm & Error & Size \\ \hline
$1$ & $15.2207$ & $20.2945$ & $1.03511$ & $8$ \\
$2$ & $1.97988$ & $1.83067$ & $0.859941$ & $23$ \\
$3$ & $0.41100$ & $0.338358$ & $0.197044$ & $69$ \\
$4$ & $0.00842346$ & $0.00815475$ & $0.00597363$ & $208$ \\
$5$ & $5.05016 \times 10^{-8}$ & $5.75969 \times 10^{-8}$ & $1.64059 \times 10^{-9}$ & $626$ \\
$6$ & $1.09651 \times 10^{-23}$ & $1.02510 \times 10^{-23}$ & $3.86286 \times 10^{-24}$
& $1878$ \\
$7$ & $1.12238 \times 10^{-70}$ & $1.22843 \times 10^{-70}$ & $8.59237 \times 10^{-71}$ & $5634$  \\
\hline
\end{tabular}
\end{center}

\medskip

\begin{center}
Method of order $4$
\end{center}

\medskip

\begin{center}
\begin{tabular}{||c|c|c|c|c||}
\hline
$n$ & $L_{2}$-norm & $L_{\infty}$-norm & Error & Size \\ \hline
$1$ & $7.44927$ & $9.64324$ & $1.04473$ & $10$ \\
$2$ & $1.81592$ & $2.41858$ & $1.15092$ & $41$ \\
$3$ & $0.0298892$ & $0.0249128$ & $0.0189903$ & $164$ \\
$4$ & $1.92454 \times 10^{-8}$ & $2.249128 \times 10^{-8}$ & $1.47053 \times 10^{-8}$ & $659$ \\
$5$ & $3.40769 \times 10^{-33}$ & $3.96407 \times 10^{-33}$ & $2.56817 \times 10^{-33}$ & $2637$ \\
\hline
\end{tabular}
\end{center}

\medskip

\begin{example}
A method of order $3$ for the evaluation of the quadratic integral
\begin{equation}\label{quadint}
I = \ifft \frac{dx}{ax^{2}+bx+c}, 
\end{equation}
\noindent
has been analyzed in \cite{manna-moll1}.  We refer to Example \ref{simpleex} 
for the explicit formulas of this Landen transformation, and define the 
iterates accordingly.  From there, we prove that the error term,
\begin{equation}
e_{n}:= \left( a_{n} - \tfrac{1}{2}\sqrt{4ac-b^{2}}, \, b_{n}, \, 
c_{n} - \tfrac{1}{2} \sqrt{4ac-b^{2}} \right)
\end{equation}
\noindent
satisfies $e_{n} \to 0$ as $n \to \infty$, with cubic rate:
\begin{equation}
\| e_{n+1} \| \leq C \|e_{n} \|^{3}.
\end{equation}
\noindent
The proof of convergence is elementary.  Therefore we have 
\begin{equation} 
(a_n,b_n,c_n) \rightarrow (\sqrt{ac-b^2/4},0,\sqrt{4ac-b^2/4})~, 
\end{equation}
which, in conjunction with (\ref{quadint}), implies that 
\begin{equation}  
I = \frac{2}{\sqrt{4ac-b^2}}\int_{-\infty}^{\infty} \frac{dx}{x^2+1}~,
\end{equation} 
exactly as one would have concluded by completing the square.  Unlike 
completing the square, our method extends to a general rational integral 
over the real line.
\end{example}

\section{The appearance of the $\, \agm $ in diverse contexts} 
\label{sec-diverse} 
\setcounter{equation}{0}

The (elliptic) Landen transformation 
\begin{equation}
a_{1}  \leftarrow \frac{1}{2}(a+b) \text{ and }
b_{1}  \leftarrow \sqrt{ab} \label{elliptic-1}
\end{equation}
\noindent
leaving invariant the elliptic integral
\begin{equation}
G(a,b) = \int_{0}^{\pi/2} \frac{ d \varphi}{\sqrt{a^{2} \cos^{2} \varphi + 
b^{2} \sin^{2} \varphi}}
\end{equation}
\noindent
appears in many different forms. In this last section we present a
partial list of them. 

\medskip

\noindent
\subsection{The Elliptic Landen Transformation} 
For the lattice 
$\mathbb{L} = \mathbb{Z} \oplus \omega \mathbb{Z},$ introduce the 
{\em theta-functions}
\begin{equation}
\vartheta_{3}(x,\omega) := \sum_{n = -\infty}^{\infty} z^{2n}q^{n^{2}} 
\text{ and }
\vartheta_{4}(x,\omega) := 
\sum_{n = -\infty}^{\infty} (-1)^{n} z^{2n}q^{n^{2}},
\end{equation}
\noindent
where $z= e^{\pi i x}$ and $q = e^{\pi i \omega}$. The condition 
$\imagpart{\omega} > 0$ is imposed to ensure convergence of the series. 
These functions admit 
a variety of remarkable identities. In particular, the {\em null-values} 
(those with $x=0$) satisfy
\begin{equation}
\vartheta_{4}^{2}(0,2 \omega) = \vartheta_{3}(0,\omega) \vartheta_{4}(0,\omega)
\text{ and } 
\vartheta_{3}^{2}(0,2 \omega) = \frac{1}{2} \left( 
\vartheta_{3}^{2}(0,\omega) + \vartheta_{4}^{2}(0,\omega) \right),
\nonumber
\end{equation}
\noindent
and completely characterize values of the $\agm$, leading to the earlier 
result \cite{borwein1}.  Grayson \cite{grayson1} has used the doubling of 
the period $\omega$ to derive the arithmetic-geometric mean  from the 
cubic equations describing the corresponding elliptic curves. See Chapter 3 
in \cite{mckmoll} for more information.  P. Sole et al. \cite{sole1, sole2}
have produced proofs of generalizations of these identities by lattice 
enumeration methods related to binary and ternary codes. \\

\noindent
\subsection{A time-one map} We now present a deeper and more modern 
version of a result known to Gauss:  
given a sequence of points $\{ x_{n} \}$ on a manifold $X$, decide whether
there is a differential equation 
\begin{equation}
\frac{dx}{dt} = V(x),
\label{diffeq}
\end{equation}
starting at $x_{0}$ such that $x_{n} = x(n,x_{0})$. Here $x(t,x_{0})$ is the 
unique solution to (\ref{diffeq}) satisfying $x(0,x_{0})=x_{0}$.  Denote by 
\begin{equation}
\phi_{ellip}(a,b) = \left( \frac{1}{2}(a+b), \sqrt{ab} \right)
\label{ellip-phi}
\end{equation}
\noindent
the familiar elliptic Landen transformation. Now take $a, b \in \mathbb{R}$ 
with $a>b>0$. Use the null-values of the theta functions to 
find unique values $(\tau, \rho)$ such that 
\begin{equation}
a = \rho \vartheta_{3}^{2}(0, \tau) \text{ and } 
b = \rho \vartheta_{4}^{2}(0, \tau).
\end{equation}
Finally define 
\begin{equation}
x_{ellip}(t) = (a(t),b(t)) = \rho \left( \vartheta_{3}^{2}(0,2^{t} \tau), 
                                    \vartheta_{4}^{2}(0,2^{t} \tau) \right),
\end{equation}
\noindent
with $x_{ellip}(0) = (a,b)$. The remarkable result is \cite{deift92}:

\begin{Thm}
(Deift, Li, Previato, Tomei). The map $t \to x_{ellip}(t)$ is an integrable 
Hamiltonian flow on $X$ equipped with an appropriate symplectic structure. 
The Hamiltonian is the complete elliptic integral $G(a,b)$ and the angle is
(essentially the logarithm of) the second period of the elliptic curve 
associated with $\tau$. Moreover 
\begin{equation}
x_{ellip}(k) = \phi_{ellip}^{k}(a,b). 
\end{equation}
Thus the arithmetic-geometric algorithm is the time-one map of a 
completely integrable Hamiltonian flow. 
\end{Thm}

Notice that this theorem shows that the result in question respects some 
additional structures whose invention postdates Gauss.

\medskip

A natural question is whether the map  (\ref{ratlanden1}) appears as a 
time-one map of an interesting flow.  \\

\noindent
\subsection{A quadruple sequence}  Several variations  of 
the elliptic Landen appear in the literature. Borchardt \cite{borchardt76}
considers  the four-term quadratically convergent iteration
\begin{eqnarray}
a_{n+1} & =  & \frac{a_{n}+b_{n} + c_{n}+d_{n}}{4}, \label{borchatdt} \\
b_{n+1} & = & \frac{\sqrt{a_{n}b_{n}} + \sqrt{c_{n}d_{n}} }{2}, \nonumber \\
c_{n+1} & = & \frac{\sqrt{a_{n}c_{n}} + \sqrt{b_{n}d_{n}} }{2}, \nonumber \\
d_{n+1} & = & \frac{\sqrt{a_{n}d_{n}} + \sqrt{b_{n}c_{n}} }{2}, \nonumber
\end{eqnarray}
\noindent
starting with $a_{0}=a, \, b_{0}=b, \, c_{0}=c$ and $d_{0}=d$. The 
common limit, denoted by $G(a,b,c,d)$, is given by
\begin{equation}
\frac{1}{G(a,b,c,d)} = \frac{1}{\pi^{2}} \int_{0}^{\alpha_{3}} 
\int_{\alpha_{1}}^{\alpha_{2}} \frac{(x-y) \, dx \, dy }{\sqrt{R(x)R(y)}},
\end{equation}
\noindent
where $R(x) = x(x-\alpha_{0})(x-\alpha_{1})(x-\alpha_{2})(x-\alpha_{3})$
and the numbers $\alpha_{j}$ are given by explicit formulas in terms of the
parameters $a, \, b, \, c,\, d$. Details are given by Mestre in 
\cite{mestre91}.

The initial conditions $(a,b,c,d) \in \mathbb{R}^4$ for which the iteration 
converges has some interesting invariant subsets.  When $a=b$ and $c=d$, we 
recover the $\agm$ iteration (\ref{agm-1}).  In the case that $b=c=d$, we 
have another invariant subset, linking to an iterative mean described below.

\medskip

\subsection
{Variations of $ \agm \, $ with hypergeometric limit}
\noindent
Let $N \in \mathbb{N}$. The analysis of 
\begin{equation}
a_{n+1} = \frac{a_{n}+(N-1)b_{n}}{N} \text{ and } c_{n+1} = 
\frac{a_{n}-b_{n}}{N},
\end{equation}
\noindent
with $b_{n} = (a_{n}^{N}-c_{n}^{N})^{1/N}$, is presented in \cite{borwein91}.
All the common ingredients appear there: a common limit, fast convergence, 
theta functions and sophisticated iterations for the evaluation of $\pi$. The 
common limit is denoted by $\text{AG}_{N}(a,b)$. The convergence is of order
$N$ and the limit is identified for small $N$: for $0 < k < 1$,
\begin{eqnarray}
\frac{1}{\text{AG}_{2}(1,k)} & = & {_2}F_{1}(1/2, \, 1/2; \, 1; \, 1 - k^{2})
\text{ and } \nonumber \\
\frac{1}{\text{AG}_{3}(1,k)} & = & {_2}F_{1}(1/3, \, 2/3; \, 1; \, 1 - k^{2}).
\nonumber 
\end{eqnarray}
\noindent
where 
\begin{equation}
{_{2}F_{1}}(a, \, b; \, c; \, x) = 
\sum_{k=0}^{\infty} \frac{(a)_{k} \, (b)_{k}}{(c)_{k} \, k!} x^{k} 
\end{equation}
\noindent
is the classical hypergeometric function.  There are integral 
representations of these as well which parallel 
(\ref{agmlaw}), see \cite{borw2}, Section 6.1 for details. \\

Other hypergeometric values appear from similar iterations. For example, 
\begin{equation}
a_{n+1} = \frac{a_{n}+3b_{n}}{4} \text{ and } 
b_{n+1} = \sqrt{b_{n}(a_{n}+b_{n})/2},
\label{itera-A4}
\end{equation}
\noindent
have a common limit, denoted by $A_{4}(a,b)$. It is given by
\begin{equation}
\frac{1}{A_{4}(1,k)}  =  {_2}F_{1}^{2}(1/4, \, 3/4; \, 1; \, 1 - k^{2}).
\end{equation}

To compute $\pi$ quartically, start at $a_{0}=1, \, 
b_{0} = (12 \sqrt{2}-16)^{1/4}$. Now compute $a_{n}$ from two steps of 
$\text{AG}_{2}$:
\begin{equation}
a_{n+1} = \frac{a_{n}+b_{n}}{2}, \text{ and } 
b_{n+1} = \left( \frac{a_{n}b_{n}^{3} + b_{n}a_{n}^{3}}{2} \right)^{1/4}.
\end{equation}
\noindent
Then
\begin{equation}
\pi = \lim\limits_{n \to \infty} 3a_{n+1}^{4} \left( 1 - \sum_{j=0}^{n} 
2^{j+1} (a_{j}^{4}-a_{j+1}^{4}) \right)^{-1}
\end{equation}
\noindent
with $| a_{n+1} - \pi | \leq C | a_{n} - \pi|^{4}$, for some constant 
$C > 0$. This is much better than
the partial sums of 
\begin{equation}
\pi = 4 \sum_{k=0}^{\infty} \frac{(-1)^{k}}{2k+1}.
\end{equation}

\medskip

The sequence $(a_{n}, \, b_{n})$ defined by the iteration 
\begin{equation}
a_{n+1} = \frac{a_{n}+2b_{n}}{3}, \quad 
b_{n+1} = \left( \frac{b_{n}(a_{n}^{2} + a_{n}b_{n} + b_{n}^{2})}{3} \right)
^{\tfrac{1}{3}}, 
\end{equation}
\noindent
starting at $a_{0} =1, \, b_{0} = x$ are analyzed in \cite{borwein89b}. They
have a  common limit $F(x)$ given by
\begin{equation}
\frac{1}{F(x)} = {_{2}F_{1}}\left( \tfrac{1}{3}, \tfrac{2}{3}; 1; 1 - x^{3} 
\right). 
\end{equation}

\medskip

\noindent
 \subsection{Iterations where the limit is harder to find} 
J. Borwein and P. Borwein \cite{borwein89c} studied the iteration of
\begin{equation}
(a,b) \to \left( \frac{a+ 3b}{4}, \frac{\sqrt{ab}+b}{2} \right),
\label{bor-1}
\end{equation}
and showed the existence of a common limit $B(a_{0},b_{0})$. Define $B(x) = 
B(1,x)$. The study of the iteration (\ref{bor-1}) is based on the functional
equation
\begin{equation}
B(x) = \frac{1+3x}{4} B \left( \frac{2(\sqrt{x}+x)}{1+3x} \right). 
\end{equation}
\noindent
and a parametrization of the iterates by theta functions \cite{borwein89c}.  The complete analysis of (\ref{bor-1}) starts with the purely computational
observation that 
\begin{equation}
B(x) \sim \frac{\pi^{2}}{3} \log^{-2}(x/4) \text{ as } x \to 0.
\end{equation}
\noindent
H. H. Chan, K. Chua and P. Sole \cite{chan} identified the limiting function 
as 
\begin{equation}
B(x) = \left( {_{2}F_{1}} \left( \tfrac{1}{3}, \, \tfrac{1}{6}; \, 1; \, 
27 \tfrac{x(1-x)^{2}}{(1+3x)^{3}} \right) \, \right)^{-2},
\end{equation}
\noindent
valid for $\tfrac{2}{3} < x < 1$. A similar hypergeometric expression gives 
$B(x)$ for $0 < x < \tfrac{2}{3}$. 

\medskip

\noindent
\subsection{Fast computation of elementary functions} The fast 
convergence of the 
elliptic Landen recurrence (\ref{agm-1}) to the arithmetic-geometric mean 
provides a method for numerical evaluation of the elliptic integral
$G(a,b)$. The 
same idea provides for the fast computation of elementary functions. For
example, in \cite{borwein3} we find the estimate
\begin{equation}
\left| \log x - \left( G(1,10^{-n}) - G(1,10^{-n}x) \right) \right|
< n 10^{-2(n-1)}, 
\end{equation}
\noindent
for $0<x<1$ and $n \geq 3$. 

\medskip

\noindent
\subsection{A continued fraction} The continued fraction
\begin{equation}
R_{\eta}(a,b) = \cfrac{a}{\eta+ \cfrac{b^{2}}{\eta+\cfrac{4a^{2}}{\eta+ 
\cfrac{9b^{2}}{\eta + \cdots} } } }, 
\end{equation}
\noindent
has an interesting connection to the $\agm$. In 
their study of the convergence of $R_{\eta}(a,b)$, J. Borwein, R. Crandall 
and G. Fee \cite{bcf1} established the identity
\begin{equation}
R_{\eta} \left( \frac{a+b}{2}, \sqrt{ab} \right) = 
\frac{1}{2} \left( R_{\eta}(a,b) + R_{\eta}(b,a) \right). 
\end{equation}
\noindent
This identity originates with Ramanujan; the 
similarity with $\agm \,$ is now direct. 

The continued fraction converges for positive real parameters, but  for 
$a, \, b \in \mathbb{C}$ the convergence question is quite delicate. For 
example, 
the even/odd parts of $R_{1}(1,i)$ converge to distinct limits. See 
\cite{bcf1} and \cite{bcf2} for more details.\\

\medskip

\noindent
\subsection{Elliptic Landen with complex initial conditions} The iteration of 
(\ref{agm-1}) with $a_{0}, \, b_{0} \in \mathbb{C}$ requires a choice of 
square root at each step. Let $a, \, b \in \mathbb{C}$ be non-zero and
assume $a \neq \pm b$. A square root $c$ of $ab$ is called the {\em right
choice} if 
\begin{equation}
\left| \frac{a+b}{2} - c \, \right| \leq \left| \frac{a+b}{2} + c \, \right|. 
\end{equation}
\noindent
It turns out that in order to have a limit for (\ref{agm-1}) one has to make
the right choice for all but finitely many indices $n \geq 1$. This is 
described in detail by Cox \cite{cox84}. 

\medskip

\noindent
\subsection{Elliptic Landen with $p$-adic initial conditions} Let $p$ be a 
prime and $a, \, b$ be 
non-zero $p$-adic numbers. In order to guarantee that the $p$-adic series
\begin{equation}
c = a \sum_{i=0}^{\infty} \binom{\tfrac{1}{2}}{i} \left( \tfrac{b}{a} -1 
\right)^{i} 
\end{equation}
\noindent
converges, and thus defines a $p$-adic square root of $ab$, one must 
assume 
\begin{equation}
b/a \equiv 1 \bmod p^{\alpha}, 
\end{equation}
\noindent
where $\alpha = 3$ for $p=2$
and $1$ otherwise. The corresponding sequence defined by (\ref{agm-1}) 
converges for $p \neq 2$ to a common limit: the $p$-adic AGM. In the case 
$p=2$ one must assume that the initial conditions satisfy $b/a \equiv 1 
\bmod 16$. In the case $b/a \equiv 1 \bmod 8$ but not $1$ modulo $16$, the 
corresponding sequence $(a_{n},b_{n})$ does not converege, but the sequence
of so-called {\em absolute invariants}
\begin{equation}
j_{n} = \frac{2^{8} (a_{n}^{4} - a_{n}^{2} b_{n}^{2} + b_{n}^{4})^{3}}
{a_{n}^{4} b_{n}^{4} (a_{n}^{2}-b_{n}^{2})^{2}}
\end{equation}
\noindent
converges to a $2$-adic integer. Information about these issues can be 
found in \cite{henniart1}. D. Kohel \cite{kohel} has proposed a generalization
of the AGM for elliptic curves over a field of characteristic
$ p \in \{ 2, \, 3, \, 5, \, 7, \, 13 \}$. Mestre \cite{mestre02} has developed 
an AGM theory for ordinary hyperelliptic curves over a field of 
characteristic $2$. This has been extended to non-hyperelliptic curves 
of genus $3$ curves by C. Ritzenhaler \cite{math.AG/0111273}. An algorithm for 
counting points for ordinary elliptic curves over finite fields of 
characteristic $p > 2$ based on the AGM is presented in R. Carls 
\cite{carls}. 

\medskip

\noindent
\subsection{Higher genus AGM} An algorithm 
analog to the AGM for abelian integrals
of genus $2$ was discussed by Richelot \cite{richelot1}, \cite{richelot2}
and Humberdt \cite{humbert}. Some details are discussed by J. Bost and J. F. 
Mestre in \cite{bost88}. The case of abelian integrals of genus $3$ is 
due to D. Lehavi \cite{math.AG/0111273} and D. Lehavi and 
C. Ritzenthaler \cite{math.AG/0403182}.

\bigskip

\noindent
{\em Gauss was correct: his numerical calculation} (\ref{once}) 
{\em has grown in
many unexpected directions}. 

\bigskip

\noindent
{\bf Acknowledgements}. The second author acknowledges the partial support of 

\noindent
$\text{NSF-DMS} 0409968$. The authors wish to thank Jon Borwein for many 
comments that lead to an improvement of the manuscript.


\bigskip


\begin{thebibliography}{10}

\bibitem{abramowitz1}
M.~Abramowitz and I.~Stegun.
\newblock {\em Handbook of {M}athematical {F}unctions with {F}ormulas, {G}raphs
  and {M}athematical {T}ables}.
\newblock Dover, New York, 1972.

\bibitem{amm1}
T.~Amdeberhan, D.~Manna, and V.~Moll.
\newblock The $2$-adic valuation of a sequence arising from a rational
  integral.
\newblock {\em Preprint}, 2007.

\bibitem{tv1}
T.~Amdeberhan and V.~Moll.
\newblock A formula for a quartic integral: a survey of old proofs and some new
  ones.
\newblock {\em The Ramanujan Journal}, 2007.

\bibitem{berndt1}
B.~Berndt and D.~Bowman.
\newblock Ramanujan's short unpublished manuscript on integrals and series
  related to {E}uler's constant.
\newblock {\em Canadian Math. Soc. Conf. Proc.}, 27:19--27, 2000.

\bibitem{borchardt76}
C.~W. Borchardt.
\newblock Uber das arithmetisch-geometrische mittel aus vier elementen.
\newblock {\em Berl. Monatsber.}, 53:611--621, 1876.

\bibitem{bomouni1}
G.~Boros and V.~Moll.
\newblock A criterion for unimodality.
\newblock {\em Elec. Jour. Comb.}, 6:1--6, 1999.

\bibitem{bomohyper}
G.~Boros and V.~Moll.
\newblock An integral hidden in {G}radshteyn and {R}yzhik.
\newblock {\em Jour. Comp. Applied Math.}, 106:361--368, 1999.

\bibitem{boros2}
G.~Boros and V.~Moll.
\newblock A rational {L}anden transformation. {T}he case of degree $6$.
\newblock In Knopp G. Mendoza E.T.~Quinto E.~L. Grinberg S. Berhanu~M, editor,
  {\em Contemporay Mathematics. Analysis, {G}eometry, {N}umber {T}heory: {T}he
  {M}athematics of {L}eon {E}hrenpreis}, volume 251, pages 83--89. American
  Mathematical Society, 2000.

\bibitem{bomoram}
G.~Boros and V.~Moll.
\newblock The double square root, {J}acobi polynomials and {R}amanujan's master
  theorem.
\newblock {\em Jour. Comp. Applied Math.}, 130:337--344, 2001.

\bibitem{boros1}
G.~Boros and V.~Moll.
\newblock Landen transformation and the integration of rational functions.
\newblock {\em Math. Comp.}, 71:649--668, 2001.

\bibitem{bomosha}
G.~Boros, V.~Moll, and J.~Shallit.
\newblock The $2$-adic valuation of the coefficients of a polynomial.
\newblock {\em Scientia}, 7:37--50, 2001.

\bibitem{borw1}
J.~M. Borwein and D.~H. Bailey.
\newblock {\em Mathematics by Experiment: Plausible reasoning in the $21$-st
  century}.
\newblock A. K. Peters, 1st edition, 2003.

\bibitem{borw2}
J.~M. Borwein, D.~H. Bailey, and R.~Girgensohn.
\newblock {\em Experimentation in Mathematics: Computational Paths to
  Discovery}.
\newblock A. K. Peters, 1st edition, 2004.

\bibitem{borwein3}
J.~M. Borwein and P.~B. Borwein.
\newblock The arithmetic-geometric mean and fast computation of elementary
  functions.
\newblock {\em SIAM Review}, 26:351--366, 1984.

\bibitem{borwein1}
J.~M. Borwein and P.~B. Borwein.
\newblock {\em Pi and the AGM- A study in analytic number theory and
  computational complexity}.
\newblock Wiley, New York, 1st edition, 1987.

\bibitem{borwein89c}
J.~M. Borwein and P.~B. Borwein.
\newblock On the mean iteration $(a,b) \leftarrow ((a+3b)/4, (
  \sqrt{ab}+b)/2)$.
\newblock {\em Math. Comp.}, 53:311--326, 1989.

\bibitem{borwein89b}
J.~M. Borwein and P.~B. Borwein.
\newblock A remarkable cubic mean iteration.
\newblock In St. Ruscheweyh, E.B. Saff, L.~C. Salinas, and R.~S. Varga,
  editors, {\em Computational {M}ethods and {F}unction {T}heory}, pages 27--31.
  Lectures Notes in Mathematics 1435, Springer-Verlag, 1990.

\bibitem{borwein91}
J.~M. Borwein and P.~B. Borwein.
\newblock A cubic counterpart of {J}acobi's identity and the {AGM}.
\newblock {\em Trans. Amer. Math. Soc.}, 323:691--701, 1991.

\bibitem{bcf1}
J.~M. Borwein, R.~Crandall, and G.~Fee.
\newblock On the {R}amanujan {AGM} fraction, {I}: the real parameter case.
\newblock {\em Experimental Mathematics}, 13:275--285, 2004.

\bibitem{bcf2}
J.~M. Borwein, R.~Crandall, and G.~Fee.
\newblock On the {R}amanujan {AGM} fraction, {II}: the complex-parameter case.
\newblock {\em Experimental Mathematics}, 13:287--295, 2004.

\bibitem{bost88}
J.~B. Bost and J.~F. Mestre.
\newblock Moyenne arithmetico-geometrique et periodes des courbes de genre $1$
  et $2$.
\newblock {\em Gazette des Mathematiciens, Soc. de Mathematique de France},
  38:36--64, 1988.

\bibitem{bromwich}
T.~J. Bromwich.
\newblock {\em An {I}ntroduction to the {T}heory of {I}nfinite {S}eries}.
\newblock Mac{M}illan, {L}ondon, 2nd edition, 1926.

\bibitem{carls}
R.~Carls.
\newblock {\em A generalized arithmetic geometric mean}.
\newblock PhD thesis, Gronigen, 2004.

\bibitem{caviness1}
B.~Caviness and J.~R. Johnson.
\newblock {\em Quantifier elimination and cylindrical algebraic decomposition}.
\newblock Texts and Monographs in Symbolic Computation. Springer-Verlag, 1st
  edition, 1998.

\bibitem{marc-moll}
M.~Chamberland and V.~Moll.
\newblock Dynamics of the degree six {L}anden transformation.
\newblock {\em Discrete and {D}ynamical {S}ystems}, 15:905--919, 2006.

\bibitem{collins2}
G.~E. Collins.
\newblock Quantifier elimination for the elementary theory of real closed
  fields by cylindrical algebraic decomposition.
\newblock {\em Lecture Notes in Computer Science}, 33:134--183, 1975.

\bibitem{cox84}
D.~Cox.
\newblock The arithmetic-geometric mean of {G}auss.
\newblock {\em L'Enseigement {M}athematique}, 30:275--330, 1984.

\bibitem{deift92}
P.~Deift.
\newblock Continuous versions of some discrete maps or what goes on when the
  lights go out.
\newblock {\em Jour. d'Analyse Math.}, 58:121--133, 1992.

\bibitem{gauss1}
K.~F. Gauss.
\newblock Arithmetisch {G}eometrisches {M}ittel.
\newblock {\em Werke}, 3:361--432, 1799.

\bibitem{grayson1}
D.~Grayson.
\newblock The arithogeometric mean.
\newblock {\em Arch. Math}, 52:507--512, 1989.

\bibitem{chan}
Kok Seng~Chua Heng Huat~Chan and P.~Sole.
\newblock Quadratic iterations to $\pi$ associated with elliptic functions to
  the cubic and septic base.
\newblock {\em Trans. Amer. Math. Soc.}, 355:1505--1520, 2002.

\bibitem{henniart1}
G.~Henniart and J.~F. Mestre.
\newblock Moyenne arithmetico-geometrique $p$-adique.
\newblock {\em C. R. Acad. Sci. Paris}, 308:391--395, 1989.

\bibitem{hubbard1}
J.~Hubbard and V.~Moll.
\newblock A geometric view of rational {L}anden transformation.
\newblock {\em Bull. London Math. Soc.}, 35:293--301, 2003.

\bibitem{humbert}
G.~Humbert.
\newblock Sur la transformation ordinaire des fonctions abeliannes.
\newblock {\em Journal des Mathematiques}, 7, 1901.

\bibitem{kauers-paule}
M.~Kauers and P.~Paule.
\newblock A computer proof of {M}oll's log-concavity conjecture.
\newblock {\em {P}reprint}, 2007.

\bibitem{kohel}
D.~Kohel.
\newblock The {AGM-X}$_{0}(${N}) {H}eegner point lifting algorithm and elliptic
  curve point counting.
\newblock In {\em Proceedings of ASI-ACRYPT'03}, volume 2894 of {\em Lecture
  Notes in Computer Science}, pages 124--136. Springer-Verlag, 2003.

\bibitem{landen1}
J.~Landen.
\newblock An investigation of a general theorem for finding the length of any
  arc of any conic hyperbola, by means of two elliptic arcs, with some other
  new and useful theorems deduced therefrom.
\newblock {\em Philos. {T}rans. {R}oyal {S}oc. {L}ondon}, 65:283--289, 1775.

\bibitem{math.AG/0111273}
D.~Lehavi.
\newblock {An explicit formula for the genus 3 AGM: arXiV:math.AG/011273}.

\bibitem{math.AG/0403182}
D.~Lehavi and C.~Ritzenhaler.
\newblock {Formulas for the arithmetic geometric mean of curves of genus 3:
  arXiV:math.AG/0403182}.

\bibitem{little}
J.~Little.
\newblock On the zeroes of two families of polynomials arising from certain
  rational integrals.
\newblock {\em Rocky Mountain Journal}, 35:1205--1216, 2005.

\bibitem{manna-moll2}
D.~Manna and V.~Moll.
\newblock Rational {L}anden transformations on $\mathbb{R}$.
\newblock {\em Math. Comp.}, 2007.

\bibitem{manna-moll1}
D.~Manna and V.~Moll.
\newblock A simple example of a new class of {L}anden transformations.
\newblock {\em Amer. Math. Monthly}, 114:232--241, 2007.

\bibitem{mckmoll}
H.~McKean and V.~Moll.
\newblock {\em Elliptic {C}urves: {F}unction {T}heory, {G}eometry,
  {A}rithmetic}.
\newblock Cambridge University Press, New York, 1997.

\bibitem{mestre91}
J.~F. Mestre.
\newblock Moyenne de {B}orchardt et integrales elliptiques.
\newblock {\em C. R. Acad. Sci. Paris}, 313:273--276, 1991.

\bibitem{mestre02}
J.~F. Mestre.
\newblock Lettre addresse a {G}audry en {H}arley.
\newblock 2000.
\newblock \url{http://www.math.jussieu.fr/~mestre}.

\bibitem{newman1}
D.~J. Newman.
\newblock A simplified version of the fast algorithm of {B}rent and {S}alamin.
\newblock {\em Math. Comp.}, 44:207--210, 1985.

\bibitem{richelot1}
F.~Richelot.
\newblock Essai sur une methode generale pour determiner la valeur des
  integrales ultra-elliptiques, fondee sur des transformations remarquables de
  ces transcendentes.
\newblock {\em C. R. Acad. Sci. Paris}, 2:622--627, 1836.

\bibitem{richelot2}
F.~Richelot.
\newblock De transformatione integralium {A}belianorum primi ordinis
  commentation.
\newblock {\em Jour. fur die reine und angew. {M}ath.}, 16:221--341, 1837.

\bibitem{sole1}
P.~Sole.
\newblock ${D}_{4}, \, {E}_{6}, \, {E}_{8}$ and the {AGM}.
\newblock {\em Springer Lecture Notes in Computer Science}, 948:448--455, 1995.

\bibitem{sole2}
P.~Sole and P.~Loyer.
\newblock ${U}_{n}$ lattices, construction {B}, and the {AGM} iterations.
\newblock {\em Europ. J. Combinatorics}, 19:227--236, 1998.

\bibitem{springer}
G.~Springer.
\newblock {\em Introduction to {R}iemann {S}urfaces}.
\newblock American {M}athematical {S}ociety, 2nd edition, 2002.

\bibitem{watson1933}
G.~N. Watson.
\newblock The {M}arquis and the {L}and-agent.
\newblock {\em Math. Gazette}, 17:5--17, 1933.

\bibitem{wilf1}
H.~S. Wilf.
\newblock {\em generatingfunctionology}.
\newblock Academic Press, 1st edition, 1990.

\end{thebibliography}
\end{document}